\def\spacingset#1{\renewcommand{\baselinestretch}{#1}\small\normalsize}

\documentstyle[12pt]{article}
\input mssymb.tex

\begin{document}
\newcommand{\nc}{\newcommand}
\nc{\rnc}{\renewcommand}
\nc{\nt}{\newtheorem}
\nc{\be}{\begin}
\nc{\erf}[1]{$\ (\ref{#1}) $}
\nc{\rf}[1]{$\ \ref{#1} $}
\nc{\lb}[1]{\mbox {$\label{#1}$} }
\nc{\hr}{\hrulefill}
\nc{\noi}{\noindent}

\nc{\eq}{\begin{equation}}
\nc{\en}{\end{equation}}
\nc{\eqa}{\begin{eqnarray}}
\nc{\ena}{\end{eqnarray}}

\nc{\ra}{\rightarrow}
\nc{\la}{\leftarrow}
\nc{\da}{\downarrow}
\nc{\ua}{\uparrow}
\nc{\Ra}{\Rightarrow}
\nc{\La}{\Leftarrow}
\nc{\Da}{\Downarrow}
\nc{\Ua}{\Uparrow}

\nc{\uda}{\updownarrow}
\nc{\Uda}{\Updownarrow}
\nc{\lra}{\longrightarrow}
\nc{\lla}{\longleftarrow}
\nc{\llra}{\longleftrightarrow}
\nc{\Lra}{\Longrightarrow}
\nc{\Lla}{\Longleftarrow}
\nc{\Llra}{\Longleftrightarrow}

\nc{\mt}{\mapsto}
\nc{\lmt}{\longmapsto}
\nc{\lt}{\leadsto}
\nc{\hla}{\hookleftarrow}
\nc{\hra}{\hookrightarrow}
\nc{\lgl}{\langle}
\nc{\rgl}{\rangle}

\nc{\stla}{\stackrel{d}{\la}}
\nc{\pard}{\partial \da}
\nc{\gdot}{\circle*{0.5}}


\rnc{\baselinestretch}{1.2}      
\nc{\bl}{\vspace{1ex}}           
\rnc{\theequation}{\arabic{section}.\arabic{equation}}  
\newcounter{xs}
\newcounter{ys}
\newcounter{os}
\nt{exc}{Exercise}
\nt{thm}{Theorem}[section]
\nt{dfn}[thm]{Definition}
\nt{pro}[thm]{Proposition}
\nt{cor}[thm]{Corollary}
\nt{con}[thm]{Conjecture}
\nt{lem}[thm]{Lemma}
\nt{rem}[thm]{Remark}
\nt{pgm}[thm]{Program}

\nc{\Poincare}{\mbox {Poincar$\acute{\rm e}$} }

\nc{\bA}{\mbox{${\bf A}$\ }}
\nc{\bB}{\mbox{${\bf B}$\ }}
\nc{\bC}{\mbox{${\bf C}$\ }}
\nc{\bD}{\mbox{${\bf D}$\ }}
\nc{\bE}{\mbox{${\bf E}$\ }}
\nc{\bF}{\mbox{${\bf F}$\ }}
\nc{\bG}{\mbox{${\bf G}$\ }}
\nc{\bH}{\mbox{${\bf H}$\ }}
\nc{\bI}{\mbox{${\bf I}$\ }}
\nc{\bJ}{\mbox{${\bf J}$\ }}
\nc{\bK}{\mbox{${\bf K}$\ }}
\nc{\bL}{\mbox{${\bf L}$\ }}
\nc{\bM}{\mbox{${\bf M}$\ }}
\nc{\bN}{\mbox{${\bf N}$\ }}
\nc{\bO}{\mbox{${\bf O}$\ }}
\nc{\bP}{\mbox{${\bf P}$\ }}
\nc{\bQ}{\mbox{${\bf Q}$\ }}
\nc{\bR}{\mbox{${\bf R}$\ }}
\nc{\bS}{\mbox{${\bf S}$\ }}
\nc{\bT}{\mbox{${\bf T}$\ }}
\nc{\bU}{\mbox{${\bf U}$\ }}
\nc{\bV}{\mbox{${\bf V}$\ }}
\nc{\bW}{\mbox{${\bf W}$\ }}
\nc{\bX}{\mbox{${\bf X}$\ }}
\nc{\bY}{\mbox{${\bf Y}$\ }}
\nc{\bZ}{\mbox{${\bf Z}$\ }}
\nc{\cA}{\mbox{${\cal A}$\ }}
\nc{\cB}{\mbox{${\cal B}$\ }}
\nc{\cC}{\mbox{${\cal C}$\ }}
\nc{\cD}{\mbox{${\cal D}$\ }}
\nc{\cE}{\mbox{${\cal E}$\ }}
\nc{\cF}{\mbox{${\cal F}$\ }}
\nc{\cG}{\mbox{${\cal G}$\ }}
\nc{\cH}{\mbox{${\cal H}$\ }}
\nc{\cI}{\mbox{${\cal I}$\ }}
\nc{\cJ}{\mbox{${\cal J}$\ }}
\nc{\cK}{\mbox{${\cal K}$\ }}
\nc{\cL}{\mbox{${\cal L}$\ }}
\nc{\cM}{\mbox{${\cal M}$\ }}
\nc{\cN}{\mbox{${\cal N}$\ }}
\nc{\cO}{\mbox{${\cal O}$\ }}
\nc{\cP}{\mbox{${\cal P}$\ }}
\nc{\cQ}{\mbox{${\cal Q}$\ }}
\nc{\cR}{\mbox{${\cal R}$\ }}
\nc{\cS}{\mbox{${\cal S}$\ }}
\nc{\cT}{\mbox{${\cal T}$\ }}
\nc{\cU}{\mbox{${\cal U}$\ }}
\nc{\cV}{\mbox{${\cal V}$\ }}
\nc{\cW}{\mbox{${\cal W}$\ }}
\nc{\cX}{\mbox{${\cal X}$\ }}
\nc{\cY}{\mbox{${\cal Y}$\ }}
\nc{\cZ}{\mbox{${\cal Z}$\ }}

\nc{\rightcross}{\searrow \hspace{-1 em} \nearrow}
\nc{\leftcross}{\swarrow \hspace{-1 em} \nwarrow}
\nc{\upcross}{\nearrow \hspace{-1 em} \nwarrow}
\nc{\downcross}{\searrow \hspace{-1 em} \swarrow}
\nc{\prop}{| \hspace{-.5 em} \times}
\nc{\wh}{\widehat}
\nc{\wt}{\widetilde}
\nc{\nonum}{\nonumber}
\nc{\nnb}{\nonumber}
 \nc{\half}{\mbox{$\frac{1}{2}$}}
\nc{\Cast}{\mbox{$C_{\frac{\infty}{2}+\ast}$}}
\nc{\Casth}{\mbox{$C_{\frac{\infty}{2}+\ast+\frac{1}{2}}$}}
\nc{\Casm}{\mbox{$C_{\frac{\infty}{2}-\ast}$}}
\nc{\Cr}{\mbox{$C_{\frac{\infty}{2}+r}$}}
\nc{\CN}{\mbox{$C_{\frac{\infty}{2}+N}$}}
\nc{\Cn}{\mbox{$C_{\frac{\infty}{2}+n}$}}
\nc{\Cmn}{\mbox{$C_{\frac{\infty}{2}-n}$}}
\nc{\Ci}{\mbox{$C_{\frac{\infty}{2}}$}}
\nc{\Hast}{\mbox{$H_{\frac{\infty}{2}+\ast}$}}
\nc{\Hasth}{\mbox{$H_{\frac{\infty}{2}+\ast+\frac{1}{2}}$}}
\nc{\Hasm}{\mbox{$H_{\frac{\infty}{2}-\ast}$}}
\nc{\Hr}{\mbox{$H_{\frac{\infty}{2}+r}$}}
\nc{\Hn}{\mbox{$H_{\frac{\infty}{2}+n}$}}
\nc{\Hmn}{\mbox{$H_{\frac{\infty}{2}-n}$}}
\nc{\HN}{\mbox{$H_{\frac{\infty}{2}+N}$}}
\nc{\HmN}{\mbox{$H_{\frac{\infty}{2}-N}$}}
\nc{\Hi}{\mbox{$H_{\frac{\infty}{2}}$}}
\nc{\Ogast}{\mbox{$\Omega_{\frac{\infty}{2}+\ast}$}}
\nc{\Ogi}{\mbox{$\Omega_{\frac{\infty}{2}}$}}
\nc{\Wedast}{\bigwedge_{\frac{\infty}{2}+\ast}}

\nc{\Fen}{\mbox{$F_{\xi,\eta}$}}
\nc{\Femn}{\mbox{$F_{\xi,-\eta}$}}
\nc{\Fp}{\mbox{$F_{0,p}$}}
\nc{\Fenp}{\mbox{$F_{\xi',\eta'}$}}
\nc{\Fuv}{\mbox{$F_{\mu,\nu}$}}
\nc{\Fuvp}{\mbox{$F_{\mu',\nu'}$}}
\nc{\cpq}{\mbox{$c_{p,q}$}}
\nc{\Drs}{\mbox{$\Delta_{r,s}$}}
\nc{\spq}{\mbox{$\sqrt{2pq}$}}
\nc{\Mcd}{\mbox{$M(c,\Delta)$}}
\nc{\Lcd}{\mbox{$L(c,\Delta)$}}
\nc{\Wxv}{\mbox{$W_{\chi,\nu}$}}
\nc{\vxv}{\mbox{$v_{\chi,\nu}$}}
\nc{\dd}{\mbox{$\widetilde{D}$}}

\nc{\diff}{\mbox{$\frac{d}{dz}$}}
\nc{\Lder}{\mbox{$L_{-1}$}}
\nc{\bone}{\mbox{${\bf 1}$}}
\nc{\px}{\mbox{${\partial_x}$}}
\nc{\py}{\mbox{${\partial_y}$}}

\setcounter{equation}{0}

\vspace{1in}
\begin{center}
{{\LARGE\bf
Moonshine Cohomology
 }}
\end{center}
\addtocounter{footnote}{0}
\vspace{1ex}
\begin{center}
Bong H. Lian and Gregg J. Zuckerman
\end{center}
\addtocounter{footnote}{0}
\footnotetext{B.H.L. is supported by grant DE-FG02-88-ER-25065.
G.J.Z. is supported by
NSF Grant DMS-9307086 and DOE Grant DE-FG02-92-ER-25121.}

\vspace{1ex}
\begin{quote}
{\footnotesize
ABSTRACT. We construct a new cohomology functor from the a certain category of
{\it quantum operator algebras} to the category of
{\it Batalin-Vilkovisky algebras}.
This {\it Moonshine cohomology} has, as a group of natural automorphisms, the
Fischer-Griess
Monster finite group. We prove a general vanishing theorem for this cohomology.
For a certain commutative QOA attached to a rank two hyperbolic lattice, we
show
that the degree one cohomology is isomorphic to the so-called Lie algebra of
physical states. In the case of a rank two unimodular lattice, the degree
one cohomology gives a new construction of Borcherd's Monster Lie algebra. As
applications,
we compute the graded dimensions and signatures of this cohomology as a
hermitean Lie algebra graded by a hyperbolic lattice.
In the first half of this paper, we give as preparations an exposition of the
theory of quantum operator algebras. Some of the results here were announced
in lectures given by the first author at the Research Institute for
Mathematical
Sciences in Kyoto in September 94.}
\end{quote}
\addtocounter{footnote}{0}

\section{ Introduction}

\ \ \ \ In the 1980's, there took place the following developments, all of
which are now understood to be connected with two-dimensional quantum field
theory:

1) Frenkel and Kac gave a construction of the simply-laced simple Lie algebras
via the vertex operator representation of the corresponding affine Kac-Moody
Lie algebras.  The vertex operators were already thought of as quantum fields
that depend holomorphically on one complex variable.  These fields had their
origin in the dual resonance model, which was central to the creation of string
theory.  Frenkel later extended this vertex operator construction to obtain
some new infinite dimensional Lie algebras that contained Kac-Moody Lie
algebras of hyperbolic type.  In the process, Frenkel gave a proof of the No
Ghost Theorem in string theory.

2) Frenkel, Lepowsky and Meurman generalized the Frenkel-Kac vertex operator
construction in order to construct an infinite dimensional graded
representation of the Monster finite group.  The existence of the FLM Moonshine
module immediately explained the empirically observed connection between the
modular function $j(\tau)$ and the dimensions of irreducible representations of
the Monster.  Borcherds soon afterwards discovered that the Moonshine module
possessed the structure of what he called a vertex algebra.  The book by
Frenkel, Lepowsky and Meurman gave a complete treatment of vertex {\it
operator} algebras and the Moonshine module.

3) Belavin, Polyakov and Zamolodchikov developed the so-called operator product
expansion into a powerful tool for the study of two dimensional conformal
quantum fields.  As a consequence, BPZ determined the possible central charges
and critical dimensions of the minimal conformal field theories. This analysis
required the deep work by Kac, and later Feigin and Fuchs on the structure of
the highest weight representations of the Virasoro Lie algebra.  It was soon
recognized (see for example \cite{MS}) that there was a close relationship
between the theory of the operator product expansion and the theory of vertex
operator algebras.

4) Feigin invented the theory of semi-infinite cohomology for Virasoro and
Kac-Moody Lie algebras.  The construction of the differential in Feigin's
complex was soon reinterpreted in terms of some elementary holomorphic quantum
fields arising in conformal field theory and the BRST quantization construction
in string theory \cite{KO}\cite{Wi85}.  Moreover, Frenkel, Garland and the
second author of the current paper were able to give a new proof of the No
Ghost Theorem via the analysis of a particular semi-infinite cochain complex.
At the end of the eighties, the two authors of the current paper began a long
series of papers building on the earlier work of FGZ.

5) Koszul discovered a new relationship between graded commutative
superalgebras and graded Lie superalgebras.  Specifically, Koszul found that
under certain very general hypotheses, the failure of a second order
differential operator to be a derivation led to the existence of a graded Lie
bracket on the underlying commutative algebra.  An important special case of
the same relationship was found independently by the physicists Batalin and
Vilkovisky.  At the time, there was no perceived connection between the work of
Koszul, Batalin and Vilkovisky and the rapidly developing study of conformal
quantum fields.

In the early 1990's, Borcherds proved the Conway-Norton conjectures for the FLM
Moonshine module.  The first two developments above are fundamental for
Borcherds.  As a brief aside, Borcherds claims in his paper that semi-infinite
cohomology theory, as discussed in 4 above, can be employed to obtain alternate
constructions of the infinite dimensional Lie algebras that figure prominently
in his work.  However, he presents no details to support his claim.

The main purpose of the current paper, partly inspired by \cite{M}, is to forge
a synthesis of Borcherds work with all five of the above developments.  We
construct a new functor that we call Moonshine cohomology and which fully
justifies Borcherds claim.  This new cohomology theory is an outgrowth of the
developments sketched in 3, 4 and 5 above.  In particular, we employ a new
mathematical approach to the operator product expansion, and we work with a
recent generalization of the notion of a VOA to the notion of a commutative
quantum operator algebra.  Our paper on CQOAs \cite{LZ11} is designed to be a
companion to the current paper on Moonshine cohomology.

Degree one Moonshine cohomology provides a functor from the category of vertex
operator algebras to the category of Lie algebras that carry an action of the
Monster finite group by automorphisms.  The degree one Moonshine cohomology of
a particular VOA yields the generalized Kac-Moody Lie algebra that is the key
to Borcherds proof of the Conway-Norton conjectures.

The total Moonshine cohomology provides a functor with values in the category
of Batalin-Vilkovsky algebras that carry an action of the Monster finite group
by automorphisms.  BV algebras are odd Poisson algebras in which the graded Lie
bracket is related to the graded commutative product in the fashion first
described by Koszul and Batalin-Vilkovisky (see 5 above.)  The abstract notion
of a BV algebra, though present in the paper of Koszul, did not become well
known until the recent work of Penkava-Schwarz and Getzler.

As in any cohomology theory, the total cohomology is better behaved and more
fundamental than the cohomology of any special degree.  Moreover, the BV
structure on the total cohomology allows us to relate a complicated Lie algebra
structure to a more elementary commutative algebra structure.  We hope that
Moonshine cohomology will yield further insight into the structure of the
Moonshine module as well as into the proof of the Conway-Norton conjectures.
We also hope that our current paper unifies a number of seemingly disparate
points of view in mathematics and mathematical physics.

Here is a brief summary of the contents of this paper:

In section 2, we state the definitions of our main concepts:  quantum
operators, matrix elements, Wick products, iterated Wick products, the
infinitely many ``circle'' products, the operator product expansion, the
notions of locality and commutativity, and finally the corresponding notions of
local and commutative quantum operator algebras. Finally, we discuss some
elementary known facts about commutativity.

In section 3, we introduce what we call the Wick calculus, which deals with
operator products of the form $:t(z)u(z):v(w)$ as well as $t(z):u(w)v(w):$
under the assumption that the quantum operators $t(z)$, $u(z)$, and $v(z)$ are
{\it pairwise commutative}. Here, $:t(z)u(z):$ denotes the Wick or normal
ordered product of $t(z)$ with $u(z)$.  The Wick calculus is essential for both
computations as well as for theoretical issues, such as the explicit
construction of CQOAs.  Section 3 continues with the construction of the CQOA
$O(b, c)$, which acts in the ghost Fock space of the BRST construction.
Following this is a construction of the CQOA $O_\kappa(L)$, which arose
originally in the seminal work of BPZ \cite{BPZ}, and which acts in the state
space of any conformal field theory having central charge $\kappa$.
This section concludes with a construction of a CQOA from a Lie algebra
equipped with an symmetric invariant form.
All three examples are special, in that these algebras are spanned by Wick
products of derivatives of the generating quantum operators. In fact, we
exhibit explicit bases consisting of such products in the first two examples.

In section 4, we discuss the BRST construction in the language of what we call
conformal QOAs.  Given a conformal QOA $O$ with central charge $\kappa$,
we form the tensor product $C^*(O) = O(b, c)\otimes O$.  We then construct the
special quantum operator $J(z)$, which we call the BRST current. We also give a
simple characterization of $J(z)$.  We then recall the famous result that the
coefficient $J(0)=Res_zJ(z)$ is square-zero if and only if $\kappa = 26$.
After that we specialize to this central charge.

The first main result of section 4 is Theorem \ref{5.3}, which states that the
Wick product induces a graded commutative associative product on the cohomology
of $C^*(O)$ with respect to the derivation, $[J(0),-]$.  This theorem first
appeared in work of E. Witten \cite{Wi3}, who called the ghost number zero
subalgebra the ``ground ring of a string background''.  An approach to this
theorem via VOA theory appears in \cite{LZ9}.  The approach in the current
paper is via CQOA theory.  Continuing section 4, we develop the theory of the
ghost field, $b(z)$, and its coefficient, $b(1)$.  As a preparation, we remind
the reader of the definition of a Batalin-Vilkovisky (BV) operator and BV
algebra. The second main result of section 4 is Theorem \ref{5.4}, which states
that the operator $b(1)$ induces a BV operator acting in the BRST cohomology
algebra. Thus the cohomology becomes a BV algebra.  This theorem was inspired
by work of Witten and Zwiebach \cite{WZ}, and first appeared i!
n \cite{LZ9}, where it was derive
At the end of section 4 we state the precise connection between BRST cohomology
and
semi-infinite cohomology.

In section 5, we finally present the construction of Moonshine cohomology
${\frak{M}}^*$
as a functor from conformal QOAs to BV algebras. Our main result is Theorem
\ref{moonvanish},
which asserts that Moonshine cohomology vanishes for degrees less than zero and
greater than three. In fact, we first prove Theorem \ref{addvanish}, which
states a vanishing theorem
for semi-infinite cohomology; Theorem \ref{moonvanish} follows immediately.

We specialize in section 6 to the Moonshine cohomology of the conformal QOA
attached by FLM to a
hyperbolic lattice of rank 2. Theorem \ref{moonshine} asserts that in this
case, the degree zero and
degree three Moonshine cohomology groups are both one dimensional;
moreover, the degree one and degree two cohomology groups can both be
canonically
identified with the so-called space of physical states, whose definition dates
back to
the early days of string theory (see \cite{Frenkel}\cite{FGZ}\cite{Bor2}). We
actually
prove a more general result, Theorem \ref{general2}, which follows from a
semi-infinite cohomology
calculation, Theorem \ref{general}.
As a consequence, we are able to compute the full
Moonshine cohomology as a tri-graded linear space, and determine all
the graded dimensions of this space. As a second application, we compute the
signature of the hermitean form on degree one cohomology and show that the form
is positive definite.

At the end of section 6 we discuss some open questions about Moonshine
cohomology.
In particular, we conjecture that the natural automorphism group of the {\it
functor} ${\frak{M}}^*$
is isomorphic to the Monster finite group. This conjecture is suggested by the
known
theorem that the Monster is the full automorphism group of the Moonshine VOA
\cite{FLM}.
We hope that future study of the functor ${\frak{M}}^*$ will cast new light on
a rich amalgam of mathematics and mathematical physics.

Sections 2-4 of this paper are meant to be an exposition, and should be
accessible to the readers who are new to the theory of operator product
expansions. Many useful (in the authors' opinion) exercises are given.
Sections 5-6 are more advanced because they draw from several different
subjects. We hope that the included references will help the readers
who are interested in further details.

{\bf Acknowledgments}: B.H.L. would like to express special thanks to
Prof. M. Miyamoto and Prof. H. Yamada for their invitation to lecture
at the Research Institute for Mathematical Sciences in Kyoto in September 94.
He also thanks them for their hospitality during his stay in Kyoto, and for
their patience waiting for the final draft of this paper. We thank F. Akman
for carefully proofreading our manuscript.

\section{Quantum Operator Algebras}\lb{sec2}

Let $V$ be a \bZ doubly graded vector space $V=\oplus V^n[m]$.
The degrees of a homogeneous element $v$ in $V^n[m]$ will
be denoted by $|v|=n$, $||v||=m$ respectively.
In physical applications, $|v|$ will be the fermion or ghost number of $v$.
In conformal field theory, $||v||$ will be the conformal dimension or weight of
$v$. We say that $V$ is bounded if for each $n$, $V^n[m]=0$ for $m<<0$.

Let $z$ be
a formal variable with degrees $|z|=0$, $||z||=-1$. Then
it makes sense to speak of a {\it homogeneous} (biinfinite)
formal power series
\eq
u(z)=\sum_{n\in\bf Z}u(n)z^{-n-1}
\en
of degrees $|u(z)|$, $||u(z)||$ where
the coefficients $u(n)$ are homogeneous linear maps in $V$ with degrees
$|u(n)|=|u(z)|$, $||u(n)||=-n-1+||u(z)||$. Note then that
the terms $u(n)z^{-n-1}$ indeed have the same degrees
$|u(z)|$, $||u(z)||$ for all $n$. {\it We require that for every $v\in V$,
$u(n)v=0$ for $n>>0$.}
(If $V$ is bounded, then this requirement is
superfluous.) We call a finite sum of such
homogeneous series $u(z)$ a {\it quantum operator} on $V$,
and we denote the linear space of quantum operators as $QO(V)$.

{\it Notations: By the expression $(z-w)^n$, $n$ an integer, we usually mean
its formal power series expansion in the region $|z|>|w|$. Thus $(z-w)^{-2}$
and $(-w+z)^{-2}$ are different, as power series. When such expressions
are to be regarded as rational functions rather than formal series, we will
explicitly mention so. When $A(z)=\sum A(n)z^{-n-1}$ is a formal series with
coefficients $A(n)$ in whatever linear space, we define $Res_zA(z)=A(0)$,
$A(z)^+=\sum_{n\geq0} A(n)z^{-n-1}$, $A(z)^-=\sum_{n<0} A(n)z^{-n-1}$,
$\partial A(z)=\sum -(n+1)A(n)z^{-n-2}$. If $u(z),u'(z)$ belong to QOAs $O,O'$
respectively, we abbreviate $u(z)\otimes u'(z)$, as an element of $O\otimes
O'$, simply as $u(z)u'(z)$. When no ambiguity occurs, we denote
$|u(z)|,||u(z)||$ simply as $|u|,||u||$. The restricted dual of a graded vector
space $V$ is denoted $V^\#$. If
$A_1(z)$,$A_2(z)$,... are quantum operators,
 an arbitrary matrix element
$\lgl x, A_1(z_1)A_2(z_2)\cdots y\rgl$ with $x\in V^\#$, $y\in V$, is denoted
as $\lgl A_1(z_1)A_2(z_2)\cdots\rgl$. In the interest of clarity, we often
write signs like $(-1)^{|t||u|}$ simply as $\pm$. This convention is used only
when the sign arises from permutation of elements. When in doubt, the reader
can easily recover the correct sign from such a permutation. Given two
homogeneous linear operators, $X,Y$, we write $[X,Y]=XY-(-1)^{|X||Y|}YX$.
A similar notation applies to quantum operators when it makes sense.}

Given two quantum operators $u(z),v(z)$, we write
\eq
:u(z)v(w):=u(z)^-v(w)+(-1)^{|u||v|}v(w)u(z)^+.
\en
Because $v(n)t=u(n)t=0$ for $n>>0$, it's easy to check that if we replace $w$
by $z$, the right hand side makes sense as a quantum operator and hence defines
a nonassociative product on $QO(V)$. It is called the {\it Wick product}.
Similarly given $u_1(z), \cdots,u_n(z)$, we define $:u_1(z_1)\cdots u_n(z_n):$
inductively as $:u_1(z_1)(:u_2(z_2)\cdots u_n(z_n):):$.
\be{exc}
Show that $:u_1(z)\cdots u_n(z):$ makes sense as an element of $QO(V)$.
\end{exc}

\be{dfn}
For each integer $n$ we define
a product on $QO(V)$:
\eq
u(w)\circ_nv(w)=Res_z u(z)v(w)(z-w)^n
-(-1)^{|u||v|}Res_zv(w)u(z)(-w+z)^n.
\en
\end{dfn}

Explicitly we have:
\eq\lb{deriv}
u(z)\circ_n v(z)=\left\{\be{array}{ll}
\frac{1}{(-n-1)!}:\partial^{-n-1}u(z) \ v(z):& \mbox{if $n<0$}\\
{[(\sum_{m=0}^n
\left(\be{array}{c}
n\\
m
\end{array}\right)
u(m)(-z)^{n-m}),v(z)]} & \mbox{if $n\geq0$.}\end{array}\right.
\en
If $A$ is a  homogeneous linear operator on $V$, then it's clear that the
graded commutator $[A,-]$ is a graded derivation of each of the products
$\circ_n$.
Since $u(z)\circ_0 v(z) = [u(0),v(z)]$, we have
\be{pro}\lb{circle0}
For any $t(z),u(z),v(z)$ in $QO(V)$ and $n$ integer, we have
\[
t(z)\circ_0(u(z)\circ_n v(z)) = {[ t(z)\circ_0 u(z)]} \circ_n v(z)  \pm
u(z)\circ_n {[t(z)\circ_0 v(z)]},
\]
ie. $t(z)\circ_0$ is a derivation of every product in $QO(V)$.
\end{pro}

\be{pro}
For $u(z),v(z)$ in $QO(V)$, the following equality of
formal power series in two variables holds:
\eq\lb{ope}
u(z)v(w)=\sum_{n\geq0}u(w)\circ_n v(w) (z-w)^{-n-1}
+:u(z)v(w):.
\en
\end{pro}
Proof: We have $u(z)v(w)=[u(z)^+,v(w)]+:u(z)v(w):$. On the other hand by
inverting the second eqn. in \erf{deriv}, we  get
\eq
[u(m),v(w)]=\sum_{n=0}^m
\left(\be{array}{c}
m\\
n
\end{array}\right)
u(w)\circ_nv(w) w^{m-n}.
\en
Thus we have
\eqa
{[u(z)^+,v(w)]}&=&\sum_{m\geq n\geq0}
\left(\be{array}{c}
m\\
n
\end{array}\right)
u(w)\circ_nv(w) w^{m-n} z^{-m-1}\nnb\\
&=&\sum_{n\geq0}u(w)\circ_nv(w)\frac{1}{n!}\partial_w^n(z-w)^{-1}\nnb\\
&=&\sum_{n\geq0}u(w)\circ_n v(w) (z-w)^{-n-1}.\ \ \ \Box
\ena

In the sense of the above Proposition, $:u(z)v(w):$ is the {\it nonsingular
part} of
the {\it operator product expansion} \erf{ope}, while
$u(w)\circ_n v(w) (z-w)^{-n-1}$ is the {\it polar part}
of order $-n-1$ (see \cite{BPZ}). In physics literature, $u(w)\circ_n v(w)$ is
often written
as $\frac{1}{2\pi i}\int_C u(z)v(w)(z-w)^ndz$ where $C$ is a small circle
around $w$.
The above proposition clearly justifies this notation.
The products $\circ_n$ will become important for describing
the algebraic and analytic structures of certain algebras
of quantum operators. Thus we introduce the following mathematical definitions:

\be{dfn}\lb{2.5}
A graded subspace \cA \ of $QO(V)$ containing the identity operator
and closed with respect to
all the products $\circ_n$ is called a quantum operator algebra.
 We say that $u(z)$ is local to $v(z)$
if $u(z)\circ_nv(z)=0$ for all but finitely many
positive $n$.
A QOA \cA is called local if its elements are pairwise mutually
local.
\end{dfn}
We observe that for any element $a(z)$ of a QOA, we have
$a(z)\circ_{-2}1 = \partial a(z)$. Thus a QOA is closed with respect
to formal differentiation.

\be{pro}\lb{2.6}
Let $u(z),v(z)$ be quantum operators, and $N$ a nonnegative integer. If
$u(z)\circ_nv(z)=0$ for $n\geq N$, then
$\lgl u(z)v(w)\rgl$ represents a rational function in $|z|>|w|$ with
poles along $z=w$ of order at most $N$.
\end{pro}
Proof: By eqn \erf{ope}, we have
\eq
\lgl u(z)v(w)\rgl = \sum_{n\geq0}\lgl u(w)\circ_n v(w)\rgl
(z-w)^{-n-1} +\lgl :u(z)v(w):\rgl.
\en
It is trivial to check that $\lgl :u(z)v(w):\rgl$, $\lgl u(w)\circ_n v(w)\rgl
\in \bC{[z^{\pm1},w^{\pm1}]}$. Thus our claim follows immediately.
$\Box$

\be{lem}\lb{operational}
Let $u(z)$ be local to $v(z)$, and $\lgl u(z)v(w)\rgl$ represent the
rational function $f(z,w)$. Then for $|w|>|z-w|$,
\eq
f(z,w)=\sum_{n\in{\bf Z}} \lgl u(w)\circ_nv(w)\rgl (z-w)^{-n-1}.
\en
\end{lem}
Proof: The Laurent polynomial $\lgl :u(z)v(w):\rgl$ in the above region is
just\\
$\sum_{i\geq0} \frac{1}{i!}\lgl :(\partial^i u(w))v(w):\rgl (z-w)^i$. Now apply
eqn. \erf{deriv}. $\Box$

We note that none of the products $\circ_n$ is associative
in general. However it clearly makes sense to speak of
the left, right or two sided ideals in a QOA as well as homomorphisms of QOAs
and they are defined in an obvious way. For example, a linear map $f:O\ra O'$
is a {\it homomorphism} if $f(u(z)\circ_n v(z))=fu(z)\circ_n fv(z)$ for all
$u(z),v(z)\in O$, and $f(1)=1$. An {\it $O$-module} is a graded space $M$
equipped with a homomorphism of QOAs $g:O\ra QO(M)$.
\be{exc}
Define the notions of a left, right, and two sided ideals for QOAs.
\end{exc}

\be{dfn}\lb{2.6b}
Two quantum operators $u(z),v(z)$ are said to commute if they are mutually
local, and
$\lgl u(z)v(w)\rgl$,$\pm \lgl v(w)u(z)\rgl$
represent the same
rational function. This is equivalent (Proposition \ref{2.6}) to the following:
for some $N\geq0$,
$(z-w)^N\lgl u(z)v(w)\rgl= \pm(z-w)^N\lgl v(w)u(z)\rgl$
as Laurent polynomials. We call a QOA $O$ whose elements pairwise commute a
commutative QOA.
\end{dfn}

\be{pro}\lb{commpro}
If $u(z),v(z)$ commute, then for all $m$
\eq
[u(m),v(w)]=\sum_{n\geq0}
\left(\be{array}{c}
m\\
n
\end{array}\right)
u(w)\circ_nv(w) w^{m-n}.
\en
\end{pro}
Proof:
The case $m\geq0$ is obtained by inverting the second eqn. in \erf{deriv}.
Since $u(z)v(w)=[u(z)^+,v(w)]+:u(z)v(w):$ and
$v(w)u(z)=\mp[u(z)^-,v(w)]\pm:u(z)v(w):$, it follows from commutativity that
$\lgl [u(z)^-,v(w)]\rgl$ represents the same rational function as
$-\lgl [u(z)^+,v(w)]\rgl$ does, which is just
$-\sum_{n\geq0}\frac{u(w)\circ_nv(w)}{(z-w)^{n+1}}$. This gives
\eq
[u(z)^-,v(w)]=-\sum_{n\geq0}u(w)\circ_nv(w)(-w+z)^{-n-1}.
\en
Taking $Res_z [u(z)^-,v(w)]z^m$ for $m<0$ gives the desired result. $\Box$

The notion of commutativity here is closely related to the physicists' notion
of duality in conformal field theory\cite{MS}. Frenkel-Lepowsky-Meurman have
reformulated
the axioms of a VOA in terms of what they call rationality,
associativity and commutativity. The notion of commutativity in Definition
\ref{2.6b} is essentially the same as FLM's. This notion has also been
reformulated in the language of formal variables in \cite{dl}.

\section{Wick's calculus}

 In this section, we derive a number of
useful formulas relating various iterated products of three quantum operators.
Most of
these formulas are well-known to physicists who are familiar with the calculus
of operator product
expansions. We will also include a lemma on commutativity.

Let $t(z),u(z),v(z)$ be homogeneous quantum operators which pairwise commute.
\be{lem}(see \cite{li})\lb{lilemma}
For all $n$, $t(z)\circ_n u(z)$ and $v(z)$ commute.
\end{lem}
Proof: We include Li's proof here for completeness. For a positive integer $N$,
$(z-w)^{2N}$ is a binomial sum of terms $(z-x)^i(x-w)^{2N-i}$, $i=1,..,2N$.
So $(z-w)^{N+2N}(t(z)\circ_n u(z)) v(w)$ is a binomial sum of terms
\eq
Res_x\left((z-w)^N(z-x)^i(x-w)^{2N-i} (t(x)u(z)(x-z)^n\mp u(z)t(x)(-z+x)^n)v(w)
\right).
\en
We want to show that for large enough $N$, and for $0\leq i\leq 2N$, term by
term we have
\eqa\lb{comm}
&&(z-w)^N(z-x)^i(x-w)^{2N-i} (t(x)u(z)(x-z)^n\mp u(z)t(x)(-z+x)^n)v(w)\nnb\\
&&=\pm(z-w)^N(z-x)^i(x-w)^{2N-i} v(w)(t(x)u(z)(x-z)^n\mp u(z)t(x)(-z+x)^n).
\ena
Consider two cases: $i\geq N$ and $i<N$.
By assumption, $(z-x)^k(t(x)u(z)(x-z)^n\mp u(z)t(x)(-z+x)^n)=0$ for all large
enough $k$. So for large enough $N$, \erf{comm} holds for $i\geq N$.
Similarly for $i<N$, $(z-w)^N(x-w)^{2N-i} (t(x)u(z)(x-z)^n\mp
u(z)t(x)(-z+x)^n)v(w)$ coincides with  $\pm(z-w)^N(x-w)^{2N-i}
v(w)(t(x)u(z)(x-z)^n\mp u(z)t(x)(-z+x)^n)$. This shows that \erf{comm} holds
for each $i$. $\Box$

This lemma is useful for showing existence of commutative QOAs: it says that
given a set of pairwise commuting quantum operators, the QOA generated by the
set  is commutative. We now develop some abstract tools for studying the
structure of commutative QOAs.

Applying \erf{ope}, we have
\eqa\lb{tuvope}
&&:t(z)u(z):v(w) \nnb\\
&=&(t(z)^-u(z) \pm u(z)t(z)^+)v(w)\nnb\\
&=&t(z)^-u(z)v(w)\pm u(z)v(w)t(z)^+ \pm u(z)[t(z)^+,v(w)]\nnb\\
&=&\sum_{n\geq0}:t(z) (u(w)\circ_nv(w)):(z-w)^{-n-1}+
:t(z)u(z)v(w):\nnb\\
&&\pm\sum_{n,m\geq0}u(w)\circ_m(t(w)\circ_nv(w))(z-w)^{-n-m-2}+\nnb\\
&&\pm\sum_{n\geq0}:u(z) (t(w)\circ_nv(w)):(z-w)^{-n-1}.
\ena
Similarly,
\eqa\lb{tuvope2}
&&t(z):u(w)v(w):\nnb\\
&=&\pm[u(w)^-,t(z)]v(w)\pm u(w)^-t(z)v(w)\pm t(z)v(w)u(w)^+\nnb\\
&=&\pm\sum_{n,m\geq0}(-1)^{n+1}
 u(w)\circ_n(t(w)\circ_mv(w))(z-w)^{-n-m-2}\nnb\\
&&\pm\sum_{n\geq0}(-1)^{n+1}:u(z)\circ_n t(z)\ v(w):(z-w)^{-n-1}\nnb\\
&&\pm\sum_{n\geq0}:u(w)\ t(w)\circ_nv(w):(z-w)^{-n-1}
\pm:u(w)t(z) v(w):
\ena

\be{lem}\lb{2.6a}
The following
equalities hold in $|w|>|z-w|$:
\eqa
(i)&& \sum_{k\in {\bf
Z}}\frac{\lgl(:t(w)u(w):)\circ_kv(w)\rgl}{(z-w)^{k+1}}\nnb\\
&&=\sum_{n,m\geq0}\frac{\lgl:\partial^mt(w)\ u(w)\circ_nv(w):\rgl
\pm\lgl:\partial^m u(w)\ t(w)\circ_nv(w):\rgl}{m!(z-w)^{n-m+1}}\nnb\\
&&\pm\sum_{n,m\geq0}
\frac{\lgl u(w)\circ_n(t(w)\circ_mv(w))\rgl}{(z-w)^{n+m+2}}\nnb\\
&&+\sum_{m\geq0}\frac{\lgl:\partial^m (t(w)u(w))\ v(w):\rgl}{m!(z-w)^{-m}}\\
(ii)&&\pm\sum_{k\in {\bf Z}}\frac{\lgl
t(w)\circ_k:u(w)v(w):\rgl}{(z-w)^{k+1}}\nnb\\
&&=\sum_{n,m\geq0}(-1)^{n+1}
\frac{\lgl u(w)\circ_n(t(w)\circ_mv(w))\rgl}{(z-w)^{n+m+2}}\nnb\\
&&+\sum_{n,m\geq0}(-1)^{n+1}\frac{\lgl:\partial^m(u(w)\circ_n t(w))v(w):\rgl}
{m!(z-w)^{n-m+1}}\nnb\\
&&+\sum_{n\geq0}\frac{\lgl:u(w)\ t(w)\circ_nv(w):\rgl}{(z-w)^{n+1}}\nnb\\
&&+\sum_{m\geq0}\frac{\lgl :u(w)(\partial^m t(w)) v(w):\rgl}{m!(z-w)^{-m}}
\ena
\end{lem}
Proof:
 To prove (i), consider  matrix coefficients on both sides of  eqn.
\erf{tuvope}. By assumption of
commutativity these matrix coefficients represent  rational functions.
Expanding both sides using Lemma \ref{operational}, we get the first eqn. (i).
The eqn. (ii) is derived similarly from \erf{tuvope2}. $\Box$

By reading off coefficients of the $(z-w)^i$, we can use this lemma to
simultaneously compute all products $:t(w)u(w):\circ_k v(w)$, and
$t(w)\circ_k :u(w)v(w):$ in terms of other products among the constituents
$t(w),u(w),v(w)$. Thus it is a kind of recursion relation for the products.
In the examples below, we will see how it allows us to understand the structure
of commutative QOAs.

\be{lem}\lb{abelianQOA}
If $t(z)^\pm u(w)^\pm=(-1)^{|t||u|}u(w)^\pm t(z)^\pm$, then\\
$:t(z)u(w)v(x):=(-1)^{|t||u|}:u(w)t(z)v(x):$.
\end{lem}
Proof: Applying the definition of the Wick product (and surpressing $z,w,x$):
\eqa
&&:tuv:-(-1)^{|t||u|}:utv:\nnb\\
&&=t^-(u^-v+(-1)^{|u||v|}vu^+)+
(-1)^{|t|(|u|+|v|)}(u^-v+(-1)^{|u||v|}vu^+)t^+\nnb\\
&&-(-1)^{|t||u|}\left(u^-(t^-v+(-1)^{|t||v|}vt^+)+
(-1)^{|u|(|y|+|v|)}(t^-v+(-1)^{|t||v|}vt^+)u^+\right)\nnb\\
&&=0.\ \ \ \Box
\ena

\subsection{Examples}

Let $QO(V)^-=\{u(z)^-\ |\ u(z)\in QO(V)\}$. This space is obviously closed
under differentiation and the Wick product. It follows that the space is also
closed under all $\circ_n$, $n$ negative. Also observe that for any
$u(z),v(z)\in QO(V)$, we have $u(z)^-v(w)^-=:u(z)^- v(w)^-:$. It follows that
the products $\circ_n$, $n=0,1,...$, restricted to $QO(V)^-$, all vanish. Thus
$QO(V)^-$ is a local QOA.

Let $LO(V)$ be the algebra spanned by homogeneous linear operators on $V$. We
can regard each operator $A$ as a formal series with just the constant term.
This makes $LO(V)$ a subspace of $QO(V)$. It is obvious that every $\circ_n$
restricted to $LO(V)$ vanishes except for $n=-1$, in which case $\circ_{-1}$ is
the usual product on $LO(V)$. Thus $LO(V)$ is a very degenerate example of a
QOA. Obviously, any commutative subalgebra of $LO(V)$ is a commutative QOA.

Let $\cC$ be the Clifford algebra with the generators $b(n),c(n)$
($n\in\bZ$) and the relations \cite{FMS}\cite{KR}\cite{Akman}
\be{eqnarray}
b(n)c(m)+c(m)b(n)&=&\delta_{n,-m-1}\nnb\\
b(n)b(m)+b(m)b(n)&=&0\nnb\\
c(n)c(m)+c(m)c(n)&=&0
\end{eqnarray}
Let $\lambda$ be a fixed integer.
The algebra \cC becomes \bZ-bigraded if we define the degrees
$|b(n)|=-|c(n)|=-1$, $||b(n)||=\lambda-n-1$, $||c(n)||=-\lambda-n$.
Let $\bigwedge^*$ be the graded irreducible
$\cC^*$-module with generator \bone\ and relations
\be{equation}
b(m)\bone=c(m)\bone=0, \ \ \ m\geq0
\end{equation}

Let $b(z),c(z)$ be the quantum operators
\be{eqnarray}
b(z)&=&\sum_mb(m)z^{-m-1}\nnb\\
c(z)&=&\sum_mc(m)z^{-m-1}
\end{eqnarray}
Let $O(b,c)$ be the smallest QOA containing $b(z),c(z)$.

\be{pro}\lb{Obc}
The QOA $O(b,c)$ is commutative. It has a basis consisting of the monomials
\eq\lb{bcmono}
:\partial^{n_1}b(z)\cdots\partial^{n_i}b(z)\
\partial^{m_1}c(z)\cdots\partial^{m_j}c(z):
\en
with $n_1>...>n_i\geq0$, $m_1>...>m_j\geq0$.
\end{pro}
Proof: Computing the OPE of $b(z),c(w)$, we have
\eqa\lb{bcope}
b(z)c(w)&=&(z-w)^{-1}+:b(z)c(w):\nnb\\
c(w)b(z)&=&(w-z)^{-1}+:c(w)b(z):\nnb\\
:b(z)c(w):&=&-:c(w)b(z):.
\ena
It follows that $b(z)$ and $c(z)$ commute. Also $b(z),c(z)$ each commutes with
itself, hence they form a pairwise commuting set.
By Lemma \ref{lilemma}, they generate a commutative QOA.

If each $u_1(z),...,u_k(z)$ is
of the form  $\partial^{n}b(z)$ or $\partial^{m}c(z)$, let's
call $:u_{1}(z)\cdots u_{k}(z):$ a monomial of degree $k$.
We claim that it's proportional to some monomial \erf{bcmono} with
$n_1>...>n_i\geq0$, $m_1>...>m_j\geq0$. If $t(z), u(z)$ each is of the form
$\partial^{n}b(z)$ or $ \partial^{m}c(z)$, it is easy to check that
$t(z)^\pm u(z)^\pm=-u(z)^\pm t(z)^\pm$. It follows from Lemma \ref{abelianQOA}
that
$:t(z)u(z)v(z):=-u(z)t(z)v(z):$ for any element $v(z)\in O(b,c)$. This shows
that $:u_{1}(z)\cdots u_{k}(z):$ is equal to
$(-1)^\sigma:u_{\sigma(1)}(z)\cdots u_{\sigma(k)}(z):$ for any permutation
$\sigma$ of $1,...,k$.

Let $O'$ be the linear span of the monomials \erf{bcmono}.
We now show that $A\circ_k B\in O'$ for any $k$ and any two monomials $A,B$,
hence $O(b,c)=O'$.
We will do a double induction on the degrees of $A$ and $B$.
Case 1: let $A=t(z)$, $B=:u(z)v(z):$ with
$t(z), u(z)$ each monomial of degree 1, and
$v(z)$ of any degree.
If $v(z)=1$, then by \erf{bcope} $t(w)\circ_k:u(w)v(z):\in O'$.
 By induction on the degree of $v(z)$ and applying
Lemma \ref{2.6a}(ii), we see that $t(w)\circ_k:u(w)v(w):\in O'$. This shows
that $A\circ_kB\in O'$ for $A$ of degree 1, $B$ of any degree.
Now case 2: suppose $A=:t(z)u(z):$, $B=v(z)$, where $t(z)$ is of degree 1 and
$u(z),v(z)$ of any degree. By induction on the degree of $u(z)$, it's clear
from Lemma \ref{2.6a}(i) that this case reduces to case 1.

Finally we must show that the monomials \erf{bcmono} are linearly independent.
 Define a map $O(b,c)\ra \bigwedge$ by $u(z)\mapsto u(-1)\bone$. We see that
this map gives a 1-1 correspondence between the set of monomials \erf{bcmono}
and a basis of $\bigwedge$. This completes the proof. $\Box$.

\be{exc}
Let $j(z)=:c(z)b(z):$. Show that $j(z)j(w)=(z-w)^{-2}+:j(z)j(w):$ by direct
computation (Hint: Eqn. \erf{tuvope} is a good guide.). Use Lemma \ref{commpro}
to conclude that $[j(m),j(n)]=n\delta_{n+m,0}$.
Construct a canonical basis of $O(j)$, the QOA generated by $j(z)$ in $O(b,c)$.
\end{exc}

Let $M(\kappa,0)$ be the Verma module of the Virasoro algebra with highest
weight $(\kappa,0)$ and vacuum vector $v_0$. Let $M(\kappa)$ be the quotient of
$M(\kappa,0)$ by the submodule generated by $L_{-1}v_0$. Let $O_\kappa(L)$ be
the QOA generated by $L(z)=\sum L_nz^{-n-2}$ in $QO(M(\kappa))$.
\be{pro}\lb{Okappa}(see \cite{BPZ}\cite{FZ})
The QOA $O_\kappa(L)$ is commutative. It has a basis consisting of monomials
\eq\lb{Lmono}
:\partial^{n_1}L(z)\cdots\partial^{n_i}L(z):
\en
with $n_1\geq...\geq n_i\geq0$.
\end{pro}
Proof: A direct computation gives
\eqa\lb{LOPE}
{[L(z)^+,L(w)]}&=&\frac{\kappa}{2}(z-w)^{-4}+2L(w)(z-w)^{-2}+\partial
L(w)(z-w)^{-1}\nnb\\
{[L(z)^-,L(w)]}&=&-\frac{\kappa}{2}(w-z)^{-4}-2L(w)(w-z)^{-2}+\partial
L(w)(w-z)^{-1}.
\ena
But we also have $L(z)L(w)=[L(z)^+,L(w)]+:L(z)L(w):$, and
$L(w)L(z)=-[L(z)^-,L(w)]+:L(z)L(w):$. Combining these with
\erf{LOPE}, it is obvious that $\lgl L(z)L(w)\rgl$ and $\lgl L(w)L(z)\rgl$
represent the same rational function.
Thus $L(z)$ commutes with itself as a quantum operator. By Lemma \ref{lilemma},
$O_\kappa(L)$ is commutative.

Let $O'$ be the linear span of the monomials \erf{Lmono} with $n_1,..,n_i\geq0$
{\it unrestricted}. To show $O'$ is closed under all the products (hence
$O_\kappa(L)=O'$), we apply induction and Lemma
\ref{2.6a} as in the case of $O(b,c)$ above. We now show that
we can restrict to those monomials \erf{Lmono} with $n_1\geq,...\geq n_i\geq0$,
and that the resulting monomials form
a basis. First by direct computation, we see (see Lemma 4.2 of \cite{LZ11})
that $O_\kappa(L)$ is a $Vir$-module defined by the action ($L(n)=L_{n-1}$)
\eq
L(n)\cdot u(z)=L(z)\circ_n u(z).
\en
Since $O_\kappa(L)$ is spanned by the monomials \erf{Lmono}, and because
$L(-n-1)\cdot u(z)=\frac{1}{n!}:\partial^{n} L(z)u(z):$ for $n\geq0$, it
follows that the module is
cyclic. Thus we have a unique onto map of $Vir$-modules $M(\kappa)\ra
O_\kappa(L)$ sending $v_0$
to $1$. But $M(\kappa)$ has a PBW basis consisting of $L(-n_1-1)\cdots
L(-n_i-1)v_0$,
$n_1\geq,...\geq n_i\geq0$. This shows that the monomials \erf{Lmono} with
$n_1\geq,...\geq n_i\geq0$
span $O_\kappa(L)$. Now define a map $O_\kappa(L)\ra M(\kappa)$ by $u(z)\mapsto
u(-1)v_0$. This is the inverse to the previous map, hence it must map a basis
to a basis. $\Box$

Let $({\frak{g}},B)$ be any Lie algebra with an invariant symmetric bilinear
form $B$,
possibly degenerate. Let $\hat{\frak{g}}$ be the affinization of
$({\frak{g}},B)$, ie. $\hat{\frak{g}}={\frak{g}}[t,t^{-1}]\oplus {\bf C}$ with
bracket:
\eq
[Xt^n,Yt^m]=[X,Y]t^{n+m}+n\delta_{n+m,0}B(X,Y)
\en
and $\bf C$ being central.
Let $M$ be any $t{\frak{g}}[t]$ locally finite $\hat{\frak{g}}$-module in which
$1\in\hat{\frak{g}}$ acts by the scalar $1$. Denote by $X(n)$ the operator
representing $Xt^n$, and define the {\it currents}:
\eq
X(z)=\sum_{n\in\bf Z}X(n)z^{-n-1}
\en
for $X\in{\frak{g}}$. Let $O$ be the QOA generated by all the currents in
$QO(M)$.
\be{pro}
The QOA $O$ is commutative. It is spanned by the monomials\\
$:\partial^{n_1}X_1(z)\cdots\partial^{n_i}X_i(z):$, with
$n_1,...,n_i\geq0$, $X_1,...,X_i\in\frak{g}$.
\end{pro}
Proof: For $X,Y\in\frak{g}$, we have
\eqa\lb{GOPE}
{[X(z)^+,Y(w)]}&=&B(X,Y)(z-w)^{-2}+[X,Y](w)(z-w)^{-1}\nnb\\
{[X(z)^-,Y(w)]}&=&-B(X,Y)(w-z)^{-2}+[X,Y](w)(w-z)^{-1}.
\ena
It follows that the currents are pairwise commuting quantum operators, and
hence $O$ is commutative.

To prove the second statement, it's enough to show that for any integer $n$ and
any two monomials $A,A'$ above, $A\circ_nA'$ is a linear sum of those
monomials.
This follows by induction on the degrees of $A,A'$ (ie. the number $X$ occuring
in $A,A'$) and by applying Lemma \ref{2.6a}. $\Box$

There is a vast literature closely related to the study of the algebra $O$
above. For a small sample, see for example
\cite{KZ}\cite{W4}\cite{TUY}\cite{FZ}\cite{Lian}\cite{li}.

\section{BRST cohomology algebras}

\be{dfn}
A conformal QOA with central charge $\kappa$ is a pair $(O,f)$, where $O$ is a
commutative QOA equipped with a homomorphism $f:O_\kappa(L)\ra O$ such that for
every homogeneous $u(z)\in O$,
\eq
fL(z)u(w) =\cdots + ||u||u(w)(z-w)^{-2} + \partial u(w) (z-w)^{-1}+:fL(z)u(w):
\en
where ``$\cdots$'' denotes the higher order polar terms. In other words,
$fL(z)\circ_1 u(z)=||u||u(w)$ and $fL(z)\circ_0 u(z)=\partial u(z)$.  For
simplicity, we sometimes write $f:O_\kappa(L)\ra O$, or simply $O_\kappa(L)\ra
O$, to denote a conformal
QOA. A homomorphism $(O,f)\stackrel{h}{\ra}(O',f')$ of conformal QOAs is a
homomorphism of QOAs $h:O\ra O'$ such that $h\circ f=f'$.
\end{dfn}

Recall that $M$ is a {\it positive energy $Vir$-module} of central charge
$\kappa$ if for every $v\in M$, we have $L_n\cdot v=0$ for $n>>0$,
and $L_0$ acts diagonalizably.
It's easy to show that $M$ is a positive energy $Vir$-module iff there's a
quantum operator $X(z)\in QO(M)$
which commutes with itself and has the OPE
\eq\lb{virope}
X(z)X(w) =\frac{\kappa}{2}(z-w)^{-4}+2X(w)(z-w)^{-2} + \partial X(w)
(z-w)^{-1}+:X(z)X(w):.
\en
\be{lem}\lb{virmodule}
Let $X(z)\in QO(M)$ define a positive energy $Vir$-module.
 Then every subalgebra of $QO(M)$ containing $X(z)$
is naturally a positive energy $Vir$-module defined by $L_n\cdot
u(z)=X(z)\circ_{n+1} u(z)$.
\end{lem}
\be{lem}\lb{virS}
Let $O$ be a commutative QOA generated by a set $S\subset QO(M)$,
and let $X(z)\in O$ have OPE \erf{virope}. Suppose for all $u(z)\in S$,
\eq\lb{Xvir}
X(z)u(w) =\cdots + ||u||u(w)(z-w)^{-2} + \partial u(w) (z-w)^{-1}+:X(z)u(w):.
\en
Then there's a unique homomorphism $f:O_\kappa(L)\ra O$ such that $fL(z)=X(z)$.
Moreover $(O,f)$ is a conformal QOA with central charge $\kappa$.
In particular, every positive energy $Vir$-module $M$ has a canonical
$O_\kappa(L)$-module
structure.
\end{lem}
We refer the interested readers to section 4 of \cite{LZ11} for the complete
proofs.

Consider, as an example, $O(b,c)$. For a fixed $\lambda$, let
\eq
X(z)=(1-\lambda):\partial b(z)\ c(z)-\lambda:b(z)\partial c(z):
\en
and $S=\{ b(z),c(z)\}$. Then we have, by direct computation \cite{FMS},
\eqa\lb{Xbc}
X(z)b(w) &=&\lambda b(w)(z-w)^{-2} + \partial b(w) (z-w)^{-1}+:X(z)b(w):\nnb\\
X(z)c(w) &=&(1-\lambda) c(w)(z-w)^{-2} + \partial c(w)
(z-w)^{-1}+:X(z)c(w):\nnb\\
X(z)X(w) &=&\frac{\kappa}{2}(z-w)^{-4}+2X(w)(z-w)^{-2} + \partial X(w)
(z-w)^{-1}+:X(z)X(w):.
\ena
where $\kappa=-12\lambda^2+12\lambda-2$. It follows that we have a conformal
QOA
$f_\lambda: O_\kappa(L)\ra O(b,c)$.

Similarly $id:O_\kappa(L)\ra O_\kappa(L)$ is itself a conformal QOA. Thus by
definition, it is the initial object in the category of conformal QOAs with
central charge $\kappa$.

\be{exc}
In the previous exercise, we define $j(z)=:c(z)b(z):$ which has $||j||=1$.
Compute the generating function $\sum_n dim\ O(j)[n]q^n$.
\end{exc}
\be{exc}
 Use Lemma \ref{virS} to classify the homomorphisms $f:O_\kappa(L)\ra O(j)$.
\end{exc}
\be{exc}
Show that $O(j)$ coincides with the subalgebra $O(b,c)^0$ of all elements of
$u(z)\in O(b,c)$ with $|u|=0$. (Hint: Compute the generating function for
the graded dimensions of $O(b,c)^0$. Alternatively for the readers who know it,
you can use the so-called boson-fermion correspondence.)
\end{exc}
\be{exc}
Use the last two exercises to classify the (for each $\lambda$) the
homomorphisms $f:O_\kappa(L)\ra O(b,c)$.
\end{exc}

\subsection{The BRST construction}

It is evident that if $(O,f)$, $(O',f')$ are conformal QOAs
on the respective spaces $V$, $V'$ with
central charges $\kappa,\kappa'$, then
$(O\otimes O',f\otimes f')$ is a conformal QOA on $V\otimes V'$
with central charge $\kappa+\kappa'$. From now on we fix $\lambda=2$ which
means that $(O(b,c),f_\lambda)$ now has central charge -26. Let $(O,f)$
be any conformal QOA with central charge $\kappa$ and consider
\eq
C^*(O)= O(b,c)\otimes O
\en
where $*$ means the total first degree.
For simplicity, we write $\hat{f}=f_\lambda \otimes f$.
\be{pro}
For every  $(O,f)$, there is a unique homogeneous element $J_f(z)\in C^*(O)$
with the following properties:\\
(i) (Cartan identity) $J_f(z)b(w)=\hat{f}L(w)(z-w)^{-1}+:J_f(z)b(w):$.\\
(ii) (Universality) If $(O,f)\ra(O',f')$ is a homomorphism of conformal QOAs,
then the induced homomorphism $C^*(O)\ra C^*(O')$ sends $J_f(z)$ to
$J_{f'}(z)$.
\end{pro}
Proof: Since the category of conformal QOAs with central charge $\kappa$ has
$(O_\kappa(L),id)$ as the initial object, if we can show that there is a unique
$J_{id}$ satisfying property (i), then (ii) implies that the same holds for
every other object in that category.

Property (i) implies $|J_{id}|=1=||J_{id}||$. Let's list a basis of the degree
$(1,1)$ subspace of $C(O_\kappa(L))$ given by Propositions \rf{Obc},
\rf{Okappa}: $:c(z)L(z):$, $:b(z)c(z)\partial c(z):$, $\partial^2c(z)$. Take a
linear combination of these elements and compute its OPE with $b(z)$. Requiring
property (i), we determine the coefficients of the linear combination and get
\eq
J_{id}(z)=:c(z) L(z):+:b(z)c(z)\partial c(z):.
\en
Now given a conformal QOA $(O,f)$, the induced map $f^*:C^*(O_\kappa(L))\ra
C^*(O)$ sends $J_{id}(z)$ to $J_f(z)=:c(z) fL(z):+:b(z)c(z)\partial c(z):$.
This completes our proof. $\Box$.

It follows from property (i) that
\eq
\hat{f}L(z)=J_f(z)\circ_0b(z)=[Q_f,b(z)]
\en
where $Q_f=Res_z J_f(z)$.
\be{lem}\lb{Qsquare}\cite{KO}\cite{Fe}\cite{FGZ}
Given $f:O_\kappa(L)\ra O$, we have
\eq\lb{JJ}
J_f(w)\circ_0J_f(w)=\frac{3}{2}\partial(\partial^2 c(w)\
c(w))+\frac{\kappa-26}{12}\partial^3 c(w)\ c(w)
\en
Thus $Q_f^2=0$ iff $\kappa=26$.
\end{lem}
Proof: We'll drop the subscripts for $J_f,Q_f$ and write $fL(z)$ as $L(z)$.
Since $J(w)\circ_0 J(w)$ is the coefficient of $(z-w)^{-1}$ in the OPE
$J(z)J(w)$, we can
extract this term from the OPE. Now $J(z)J(w)$ is
the sum of 4 terms:
\eqa
&(i)& c(z) L(z)c(w)L(w)\nnb\\
&(ii)& c(z)L(z):b(w)c(w)\partial c(w):\nnb\\
&(iii)& :b(z)c(z)\partial c(z): c(w)L(w)\nnb\\
&(iv)& :b(z)c(z)\partial c(z)::b(w)c(w)\partial c(w):.
\ena
 Extracting the coefficient of $(z-w)^{-1}$ (which is done by applying Lemma
\ref{2.6a} repeatedly) in each of these 4 OPEs, we
get respectively (surpressing $w$):
\eqa
&(i)& 2\partial c\ cL+\frac{\kappa}{12}\partial^3c\ c\nnb\\
&(ii)& c\partial c\ L\nnb\\
&(iii)& c\partial c\ L\nnb\\
&(iv)& \frac{3}{2}\partial(\partial^2 c\ c)-\frac{13}{6}\partial^3 c\ c
\ena
Thus \erf{JJ} follows. Now $2Q^2=[Q,Q]=Res_w[Q,J(w)]=Res_w J(w)\circ_0 J(w)$,
which is zero iff $\kappa=26$. $\Box$

Given $f:O_{26}(L)\ra O$, $[Q_f,-]=J_f(z)\circ_0$ is a derivation of the QOA
$C^*(O)$ (Lemma \ref{circle0}). For $\kappa=26$, which we assume from now on,
$[Q_f,-]$ becomes
a differential on $C^*(O)$ and we have a cochain complex
\eq
[Q_f,-]:C^*(O)\lra C^{*+1}(O).
\en
It is called {\it the BRST complex} associated to $f:O_{26}(L)\ra O$. Its
cohomology will be denoted as $H^*(O)$. All the operations $\circ_n$ on
$C^*(O)$
descend to the cohomology. However, all but one is trivial.
\be{thm}\lb{5.3}\cite{Wi3}\cite{WZ}\cite{LZ9}
The Wick product $\circ_{-1}$ induces a graded
commutative associative product on
$H^*(O)$ with unit element represented by the identity operator. Moreover,
every cohomology class is represented by a quantum operator $u(z)$ with
$||u||=0$.
\end{thm}

\be{exc}
Check that $1$ represents the unit of the commutative algebra $H^*(O)$. Show
that for all $n\neq-1$, $\circ_n$ is homologically trivial, ie. if $u(z),v(z)$
represent two cohomology classes, then $u(z)\circ_n v(z)=[Q_f,t(z)]$ for some
$t(z)$. (Hint: Recall that $\partial A(w)=[Q_f,b(w)]\circ_0A(w)$ and
$||A||A(w)=[Q_f,b(w)]\circ_1A(z)$.)
\end{exc}

\subsection{Batalin-Vilkovisky Algebras}

Let $A^*$ be a \bZ graded commutative associative algebra.
For every $a\in A$, let $l_a$ denote the linear map on $A$  given by
the left multiplication by $a$. Recall that a (graded) derivation $d$ on $A$ is
a homogeneous linear operator such that ${[d,l_a]} - l_{da}=0$ for all $a$.
A BV operator \cite{W2}\cite{SP}\cite{GJ} $\Delta$ on $A^*$ is
a linear operator of degree -1 such that:\\
(i) $\Delta^2=0$;\\
(ii) ${[\Delta,l_a]} - l_{\Delta a}$ is a derivation on $A$ for all $a$, ie.
$\Delta$ is
a {\it second} order derivation.

A BV algebra is a pair $(A,\Delta)$ where $A$ is a graded commutative
algebra and $\Delta$ is a BV operator on $A$. The following is an elementary
but fundamental lemma:
\be{lem}\cite{Koszul}\cite{GJ}\cite{Pen}
Given a BV algebra $(A,\Delta)$, define the BV bracket $\{,\}$ on $A$ by:
\[
(-1)^{|u|}\{u,v\} = {[\Delta,l_u]}v - l_{\Delta u}v.
\]
Then $\{,\}$ is a graded Lie bracket on $A$ of degree -1, ie.
\eqa
&&\{u,v\}+(-1)^{(|u|-1)(|v|-1)}\{v,u\}=0\nnb\\
&&(-1)^{(|u|-1)(|t|-1)}\{u,\{v,t\}\}
+(-1)^{(|t|-1)(|v|-1)}\{t,\{u,v\}\}
+(-1)^{(|v|-1)(|u|-1)}\{v,\{t,u\}\}=0\nnb
\ena
\end{lem}
By property (ii) above, it follows immediately that for every
$u\in A$, $\{u,-\}$ is a derivation on $A$. Thus a BV algebra is a special kind
of an odd Poisson algebra which, in mathematics, is also known as a
Gerstenhaber
algebra \cite{Gers1}. It's important to note
 that {\it $A^1$ is canonically a Lie algebra and that each $A^p$ is an
$A^1$-module}.

\be{exc}
Let $\frak{g}$ be an Lie algebra, $\wedge^*\frak{g}$ its exterior
algebra and $\delta$ the Chevalley-Eilenberg differential on
$\wedge^*\frak{g}$.
Check that $(\wedge^*\frak{g},\delta)$ is a BV algebra. Show that
$\wedge^1\frak{g}$ is canonically isomorphic to $\frak{g}$ as a Lie algebra.
\end{exc}
\be{exc}
 More generally,
let $B$ be any commutative algebra and $f: {\frak{g}}\ra Der\ B$ be a Lie
algebra homomorphism (making $B$ a $\frak{g}$-module). Consider the
Lie algebra homology complex ${\wedge^*\frak{g}}\otimes B$. Show that the
Chevalley-Eilenberg differential is a BV operator on this complex.
\end{exc}

Given $f:O_{26}(L)\ra O$, consider  the linear operator $\Delta_f: C^*(O)\lra
C^{*-1}(O)$,
$u(z)\mapsto b(z)\circ_1 u(z)$.
\be{thm}\lb{5.4}\cite{LZ9}
The operator $\Delta_f$ descends to the cohomology $H^*(O)$. Morever,
it is a BV operator on the commutative algebra $H^*(O)$. Thus
$H^*(O)$ is naturally a BV algebra.
\end{thm}
For complete proofs of the two theorems above, see section 4 of \cite{LZ11}.
The theorems were originally proved in \cite{LZ9} in the
context of vertex operator algebras. (For related
versions of Theorem \ref{5.4}, see \cite{GJ}\cite{SP}\cite{KSV}\cite{Hu}.)

\be{pgm}\lb{4.9}
 Study $H^*(O\otimes O')$ as a bifunctor from pairs of conformal
QOAs to BV algebras. In particular, fix $O$ and study $H(O\otimes -)$ as a
functor from conformal QOAs to BV algebras. An automorphism group of $O$ acts
by natural automorphisms on the functor $H(O\otimes -)$.
\end{pgm}

\subsection{Modules}

Consider $f:O_{26}(L)\ra O$, and an $O$-module $M$ equipped with the structure
homomorphism $g:O\ra QO(M)$. The homomorphism $g$ induces
$g^*:C^*(O)\ra QO(\wedge\otimes M)$. We write $g^*J_f(z)$ as $J_{f,g}(z)$, its
residue as $Q_{f,g}$, and  the $C^*(O)$-module $\wedge\otimes M$ as $C^*(O,M)$.
By Lemma \ref{Qsquare}, $Q_{f,g}$ turns $C^*(O,M)$ into a complex whose
cohomology is denoted as $H^*(O,M)$. It turns out that $H^*(O,M)$ is a module
over the BV algebra $H^*(O)$ in a suitable sense. This will be the topic
of a future paper.

Let $N$ be a positive energy $Vir$-module of central charge 26.  By Lemma
\ref{virS},
we have a canonical homomorphism $g:O_{26}(L)\ra QO(N)$.
This makes $\wedge\otimes N$ into a $C^*(O_{26}(L))$-module with BRST
differential $Q_{id,g}$. On the other hand, $\wedge\otimes N$ is by
definition the space of semi-infinite cochains $C^{\frac{\infty}{2}+*}(Vir,N)$
of $Vir$ with coefficients in $N$ (see \cite{Fe}\cite{FGZ}). The differential
of this cochain complex is denoted by $d_N$, and its cohomology as
$H^{\frac{\infty}{2}+*}(Vir,N)$.
It's easily seen that we have $d_N=Q_{id,g}$ \cite{FGZ}. It follows that we
have
\eq
H^*(O_{26}(L),N)\cong H^{\frac{\infty}{2}+*}(Vir,N).
\en

Now given a conformal QOA $(O,f)$ and an $O$-module $M$ equipped with the
homomorphism $g:O\ra QO(M)$, we can regard $M$ as an $O_{26}(L)$-module via
$g\circ f:O_{26}(L)\ra QO(M)$.
It follows that $(C^*(O,M),Q_{f,g})=(C^*(O_{26}(L),M), Q_{id,g\circ f})=
(C^{\frac{\infty}{2}+*}(Vir,M),d_M)$ as complexes.

Recall the linear isomorphism $O(b,c)\stackrel{\sim}{\ra}\wedge$, $u(z)\mapsto
u(-1)\bone$ (see proof of Lemma \ref{Obc}). Given a conformal QOA $(O,f)$, we
have $C^*(O)=O(b,c)\otimes O\stackrel{\sim}{\ra}\wedge\otimes O$. Call this map
$k$.
 We claim that the differential induced on
$\wedge\otimes O$ via $k$ coincides with the semi-infinite differential $d_O$.
We must check that
$J_f(z)\circ_0 (u(z)\otimes v(z))\stackrel{k}{\mapsto} d_O (u(-1)\bone\otimes
v(z))$.

 Let  $a(z)=:b(z)c(z)\partial c(z):$. It acts only on $\wedge$, and hence
$a(z)\circ_0(u(z)\otimes v(z))={[a(0),u(z)]}
\otimes v(z)\stackrel{k}{\mapsto}a(0)u(-1)\bone\otimes v(z)$.
Using Lemma \ref{operational} to compute the OPE $c(z)fL(z)(u(w)\otimes v(w))$
and get
$c(z)fL(z)\circ_0(u(z)\otimes v(z))=\sum (c(z)\circ_nu(z)) \otimes
(fL(z)\circ_{-n-1}v(z))$. It's also easy to
check, using eqn.\erf{deriv}, that under $O(b,c)\ra \wedge$, we have
$c(z)\circ_n u(z)\mapsto c(n)u(-1)\bone$. It follows that
\eqa
J_f(z)\circ_0(u(z)\otimes v(z))\stackrel{k}{\mapsto}&&\sum
c(n)u(-1)\bone\otimes
L(-n-1)\cdot v(z)+a(0)u(-1)\bone\otimes v(z)\nnb\\
&&=d_O ( u(-1)\bone\otimes v(z)).
\ena
The last equality follows from the definition of $d_O$. Thus we've shown that
\be{lem}\lb{semiinf}
Given a conformal QOA $(O,f)$ and an $O$-module $M$ equipped with homomorphism
$g:O\ra QO(M)$, we have\\
(i) $(C^*(O),Q_f)\cong (C^{\frac{\infty}{2}+*}(Vir,O),d_O)$.\\
(ii)$(C^*(O,M),Q_{f,g})=(C^{\frac{\infty}{2}+*}(Vir,M),d_M)$.
\end{lem}

\section{Moonshine Cohomology}

We'll now construct a functor $\frak{M}$ using the Moonshine VOA of rank 24 and
the BRST construction above.
This will be a functor from the category of conformal QOAs of central charge 2
to the category of BV algebras (Program \ref{4.9}).
In particular, it assigns a Lie algebra ${\frak{M}}^1(O)$ to every such
conformal QOA $O$. As a special case, if
$O$ is the conformal QOA corresponding to the unimodular rank 2 hyperbolic
lattice $II_{1,1}$ (see below), then ${\frak{M}}^1(O)$ is Borcherds' Monster
Lie algebra.

Let $(V,\bone,\omega,Y(-,z))$ be the Moonshine VOA as studied by
Frenkel-Lepowsky-Meurman and Borcherds \cite{FLM2}\cite{Bor}\cite{FLM}
(see Definitions 8.10.1-8.10.18 of \cite{FLM}).
It's now known that\\
(i) The Fischer-Griess Monster finite group $F_1$ is the automorphism group of
the VOA $V$.\\
(ii) $\sum dim\ V[n]q^{n-1}=j(q)-744$ where $j(q)$ is the Dedekind-Klein
j-function.\\
(iii) $Y(\omega,z)$ defines a unitarizable $Vir$-module structure on $V$.

\be{pro} (see \cite{LZ10})\lb{voatoqoa}
Let $(U,\bone,\omega,Y(-,z))$ be VOA of rank $\kappa$. Let $O(U)$ be the linear
space of vertex operators, ie. $O(U)=\{Y(a,z)|a\in U\}$. Then $O(U)\subset
QO(U)$ is a conformal QOA with $O_\kappa(L)\ra O(U)$,
$L(z)\mapsto Y(\omega,z)$. Moreover, $O(U)$ has an $O(U)$-module structure
$O(U)\ra QO(O(U))$ defined by $u(z)\mapsto \hat{u}(z)$, $\hat{u}(n)\cdot
v(z)\stackrel{def}{=}u(z)\circ_n v(z)$.
\end{pro}
Proof:
 By Lemma \ref{operational} above and Proposition 8.10.5 of \cite{FLM}, we
have\\ $\sum\lgl Y(u,w)\circ_n Y(v,w)\rgl (z-w)^{-n-1}=\lgl
Y(Y(u,z-w)v,w)\rgl$. Thus $O(U)$ is closed under all the products $\circ_n$.
We also see that $O(U)$ has an $O(U)$-module structure as claimed.
By Proposition 8.10.3 of \cite{FLM}, $O(U)$ is commutative. By definition,
$Y(\omega,z)=\sum\omega(n)z^{-n-1}$ satisfies, for all $u$,
$Y(\omega(0)u,z)=\partial Y(u,z)$, $\omega(1)u=||u||u$,
$\omega(2)\omega=\frac{\kappa}{2}\omega$. By Lemma \ref{virS}, this
means that $O(U)$ is a conformal QOA as claimed. $\Box$

It follows immediately that the linear bijection $U\ra O(U)$,
$a\mapsto Y(a,z)$ is an isomorphism of $O(U)$-modules. It's also clear that any
automorphism of
the VOA $U$ yields an automorphism of the conformal QOA $O(U)$. In the case of
the
Moonshine VOA $V$, it follows that\\
(iv) $O(V)$ is a conformal QOA of central charge 24, in which $F_1$ acts by
automorphisms.

Thus for any conformal QOA $O$ of central charge 2, we can consider the BRST
QOA
$C^*(O(V)\otimes O)$. We denote its cohomology as ${\frak{M}}^*(O)$, which by
Theorem \ref{5.4}
is a BV algebra. We call ${\frak{M}}^*(O)$ {\it the Moonshine cohomology of
$O$}.
 If $M$ is an $O$-module, then $V\otimes M$ is naturally an $O(V)\otimes
O$-module.
It follows from the preceeding section that $C^*(O(V)\otimes O,V\otimes M)$ is
a cochain complex.
We denote its cohomology as ${\frak{M}}^*(O,M)$, which we call {\it the
Moonshine cohomology
of $(O,M)$}.
It follows from (iv) above that $F_1$ is a group of automorphisms of both
${\frak{M}}^*(O)$ and ${\frak{M}}^*(O,M)$ (Program \ref{4.9}).

\subsection{Vanishing Theorem}

A $Vir$-module is {\it tame} if it's graded dimensions are finite.
A $Vir$-module is {\it hermitean} if it's a direct sum of a tame
 positive energy modules equipped with
an invariant nondegenerate hermitean form.
A hermitean $Vir$-module is {\it unitarizable} if its hermitean form is
positive definite.
{\it Unless specified otherwise, $Vir$-modules and conformal QOAs from now on
are assumed to have the first degree $|\cdot|\equiv0$.}

\be{thm}\lb{moonvanish}
For any conformal QOA $O$ of central charge 2, and any $O$-module $M$,
${\frak{M}}^p(O)$, ${\frak{M}}^p(O,M)$ vanish for all $p\neq0,1,2,\ or\ 3$.
\end{thm}
\be{thm}\lb{addvanish}
Let $r$ be a real number with $1<r<25$.
Let $P$ and $N$ be positive energy $Vir$-modules of central charges $26-r,r $
respectively. Assume $P$ is unitarizable.
Then $H^{\frac{\infty}{2}+p}(Vir,P\otimes N)$ vanishes for all $p\neq0,1,2,\
or\ 3$.
\end{thm}
Proof: By the unitarizability of $P$, it's a direct sum of irreducible modules
$L(26-r,h)$, $h\geq0$, with suitable multiplicities.
Thus it's enough to do the case $P=L(26-r,h)$.

 Recall that $N$ is a $Vir$-module of central charge $r$, in which $L_0$ acts
diagonalizably.
Since every cohomology class in $H^{\frac{\infty}{2}+*}(Vir,L(26-r,h)\otimes
N)=0$ is represented by an element of zero weight, we may assume, without loss
of generality, that $L_0$ only has real eigenvalues in $N$. Thus any
irreducible module $L(r,k)$ occuring in the composition series of $N$ must have
real $k$.
{}From the structure of the Verma modules, $L(r,k)=M(r,k)$ unless $k=0$, and
$L(26-r,h)=M(26-r,h)$ for $h>0$.

By our reduction theorem on semi-infinite cohomology \cite{LZ2}, for $k\neq0$,
or $h>0$, we have $H^{\frac{\infty}{2}+p}(Vir,L(26-r,h)\otimes L(r,k))=0$
for $p\neq1,2$.
It's easy to verify that \\ $H^p(O_{26}(L),L(26-r,0)\otimes L(r,0))$ is zero if
$p\neq0,3$, and one dimensional if $p=0,3$. Thus if $N$ is a module of finite
length, then
$H^{\frac{\infty}{2}+p}(Vir,L(24,0)\otimes N)=0$ for $p\neq0,..,3$.
Every finitely generated positive energy $Vir$-module of
central charge $r$ has finite length. Now any module is a direct limit
of finitely generated modules, and direct limit is exact with respect to
cohomology.
If follows that $H^{\frac{\infty}{2}+p}(Vir,L(24,0)\otimes N)=0$ for
$p\neq0,..,3$.
This completes the proof. $\Box$

Proof of Theorem \ref{moonvanish}:
Specialize Theorem \ref{addvanish} to the case $r=2$, $P=O(V)$ (which is
unitarizable), $N=O$, and applying
Lemma \ref{semiinf}, we see that ${\frak{M}}^p(O)$ vanishes for all
$p\neq0,1,2,\ or\ 3$. For $N=M$, we have a similar statement for
${\frak{M}}^p(O,M)$. $\Box$

\section{Moonshine cohomology and the Monster Lie algebra}

Let $M$ be a positive energy $Vir$-module. Let $Vir_\pm$ be respectively the
subalgebras spanned by
the $L_n$, $\pm n>0$. Define the {\it physical space associated to $M$}:
\eq
{\frak{P}}(M)=M[1]^{Vir_+}/N(M)
\en
where $N(M)=(Vir_-\cdot M)\cap M[1]^{Vir_+}$. If $M$ is a hermitean module of
central charge 26, then
there are two natural linear maps
$\nu_i:{\frak{P}}(M)\ra H^{\frac{\infty}{2}+i}(Vir,M)$, $i=1,2$, given
respectively by $v\mapsto c(-1)v$, $v\mapsto c(-2)c(-1)v$ (see \cite{LZ9}
section 2.4 for details). To emphasize their dependence on $M$, we'll refer to
these maps as {\it $\nu_1,\nu_2$ for the module $M$}.
Let $O$ be a conformal QOA.
Suppose the $Vir$-module structure $f:O_\kappa(L)\ra QO(O)$ on $O$, given by
Lemma \ref{virmodule}, is hermitean. Then it makes sense to
consider the maps $\nu_1,\nu_2$ for $O$.

If $u(z)\in O^{Vir_+}$ then it's easy to show using commutativity that\\
 $u(z)\circ_0fL(z)=fL(z)\circ_0u(z)-\partial u(z)=0$. This implies that
\eq
u(z)\circ_0(fL(z)\circ_nv(z))=fL(z)\circ_n(u(z)\circ_0v(z)).
\en
For $v(z)\in O^{Vir_+}$, this shows that $u(z)\circ_0v(z)\in O^{Vir_+}$. For
$v(z)\in N(O)$, it shows
that $u(z)\circ_0v(z)\in N(O)$. Using commutativity, we show that
$u(z)\circ_0v(z)+v(z)\circ_0u(z)=\partial A(z)$
for some $A(z)$. Thus $\circ_0$ is a skew symmetric product on $O^{Vir_+}$
modulo $N(O)$, and it also
factors through $N(O)$. The fact that $u(z)\circ_0$ is a derivation of the
product $v(z)\circ_0t(z)$
says exactly that the skew symmetric operation $\circ_0$ on ${\frak{P}}(O)$
satisfies that Lie algebra
Jacobi identity. Thus ${\frak{P}}(O)$ is a Lie algebra with bracket $\circ_0$.
We'll use the convention
that $-u(z)\circ_0v(z)$ is the Lie bracket of $u(z)$ with $v(z)$.

If $O$ has central charge 26, then the maps $\nu_i$ together with Lemma
\ref{semiinf} yield
two new maps (which we also call $\nu_i$), $\nu_i:{\frak{P}}(O)\ra H^i(O)$.
The bracket in $H^1(O)$ can be written as
$\{A(z),B(z)\}=(-1)^{|A|}(b(z)\circ_0A(z))\circ_0 B(z)$.
(see \cite{LZ11} section 5 for details). Thus,
\eqa
\{\nu_1 u(z),\nu_1 v(z)\}&=&\{c(z)u(z),c(z)v(z)\}\nnb\\
&=&-u(z)\circ_0 (c(z) v(z))\nnb\\
&=&\nu_1(-u(z)\circ_0v(z)).
\ena
Thus $\nu_1$ is a Lie algebra homomorphism. Since $H^2(O)$ is canonically a
$H^1(O)$-module,
it becomes a ${\frak{P}}(O)$-module via $\nu_1$. But we also have
\eq
\{\nu_1 u(z),\nu_2 v(z)\}=\nu_2 (-u(z)\circ_0v(z)).
\en
Thus $\nu_2$ is a ${\frak{P}}(O)$-module homomorphism from the adjoint module
${\frak{P}}(O)$ to $H^2(O)$.
To emphasize their dependence on the QOA $O$, we'll refer to those
homomorphisms as {\it $\nu_1$, $\nu_2$ for the QOA $O$}.
To summarize,
\be{lem}\lb{lahomo}
Given a hermitean $Vir$-module $M$ of central charge 26, we've two linear maps
$\nu_i:{\frak{P}}(M)\ra H^{\frac{\infty}{2}+i}(Vir,M)=H^i(O_{26}(L),M)$,
$i=1,2$, given by
$u\mapsto c(-1)u$, $u\mapsto c(-2)c(-1)v$ respectively.
Given a hermitean conformal QOA $O$ of central charge 26, we've a Lie algebra
homomorphism and
a module homomorphism $\nu_i:{\frak{P}}(O)\ra H^i(O)$, $i=1,2$,
given by $u(z)\mapsto c(z)u(z)$, $u(z)\mapsto \partial c(z)\ c(z)u(z)$
respectively.
\end{lem}

Let $\Lambda$ be any rank 2 hyperbolic even integral lattice, and
$(V_\Lambda,\bone_\Lambda,
\omega_\Lambda,Y_\Lambda(-,z))$ be the canonical rank 2 VOA associated to
$\Lambda$ \cite{FLM}. The $Vir$-module $O(V_\Lambda)\cong V_\Lambda$ is a
direct sum of the so-called Fock modules, which are hermitean.
By the lemma above, we have
$\nu_i:{\frak{P}}(O(V)\otimes O(V_\Lambda))\ra H^i(O(V)\otimes O(V_\Lambda))
={\frak{M}}^i (O(V_\Lambda))$, $i=1,2$.
\be{thm}\lb{moonshine}
The homomorphisms $\nu_1$, $\nu_2$
for the QOA $O(V)\otimes O(V_\Lambda)$ are isomorphisms, and we've
\eq
{\frak{M}}^p (O(V_\Lambda))=\left\{\be{array}{ll}
\bC 1& \mbox{if $p=0$}\\
\nu_1 {\frak{P}}(O(V)\otimes O(V_\Lambda)) & \mbox{if $p=1$.}\\
\nu_2 {\frak{P}}(O(V)\otimes O(V_\Lambda)) & \mbox{if $p=2$.}\\
\bC \partial^2c(z)\partial c(z) c(z) & \mbox{if $p=3$}\\
0& \mbox{otherwise.}\end{array}\right.
\en
\end{thm}
\be{cor}
Let $\Lambda$ be the unimodular lattice $II_{1,1}$. Then ${\frak{M}}^1
(O(V_\Lambda))$ is canonically
isomorphic to the Monster Lie algebra, and ${\frak{M}}^2 (O(V_\Lambda))$ to the
adjoint module.
\end{cor}
Proof: By definition \cite{Bor2}, the Monster Lie algebra has as its underlying
space
${\frak{P}}(V\otimes V_\Lambda)$, and its bracket $[u,v]=-Res_z Y(u,z)v$. Now
${\frak{P}}(O(V)\otimes O(V_\Lambda))\cong {\frak{P}}(V\otimes V_\Lambda)$
follows
from Proposition \ref{voatoqoa}. $\Box$

The rest of this paper is devoted to proving and generalizing the theorem
above.

\subsection{Hyperbolic lattices}

Let $\Lambda$ be a rank $r\leq26$ even integral Lorentzian lattice, and
$(V_\Lambda,\bone_\Lambda,\omega_\Lambda,Y_\Lambda(-,z))$ the canonical rank
$r$ VOA associated to $\Lambda$. This VOA has rank $r$. Let $O$ be any
conformal QOA
of central charge $26-r$ such that it's unitarizable as a $Vir$-module.
\be{thm}\lb{general2}
Under the above assumptions the homomorphisms
 $\nu_1,\nu_2$ for the QOA $O\otimes O(V_\Lambda)$ are isomorphisms, and
\eq
H^p(O\otimes O(V_\Lambda))=\left\{\be{array}{ll}
Hom_{Vir}(L(26-r,0),O) & \mbox{if $p=0$}\\
\nu_1 {\frak{P}}(O\otimes O(V_\Lambda)) & \mbox{if $p=1$.}\\
\nu_2 {\frak{P}}(O\otimes O(V_\Lambda)) & \mbox{if $p=2$.}\\
\partial^2c(z)\partial c(z) c(z)Hom_{Vir}(L(26-r,0),O) & \mbox{if $p=3$}\\
0& \mbox{otherwise.}\end{array}\right.
\en
\end{thm}
Theorem \ref{moonshine} is clearly an immediate consequence when $r=2$ and
$O=O(V)$.
By Lemma \ref{semiinf}, Theorem \ref{general2} is equivalent to

\be{thm}\lb{general}
Under the above assumptions the linear maps
 $\nu_1,\nu_2$ for the module $O\otimes V_\Lambda$ are bijections, and
\eq
H^{\frac{\infty}{2}+p}(Vir,O\otimes V_\Lambda)=\left\{\be{array}{ll}
Hom_{Vir}(L(26-r,0),O) & \mbox{if $p=0$}\\
\nu_1 {\frak{P}}(O\otimes V_\Lambda) & \mbox{if $p=1$.}\\
\nu_2 {\frak{P}}(O\otimes V_\Lambda) & \mbox{if $p=2$.}\\
c(-3)c(-2)c(-1)Hom_{Vir}(L(26-r,0),O) & \mbox{if $p=3$}\\
0& \mbox{otherwise.}\end{array}\right.
\en
\end{thm}
We now proceed to prove this. Let $\bR^{k,l}$ be the standard pseudo-euclidean
space of signature
$(k,l)$. The inner product is written as $\alpha\cdot\alpha$.
Given $\alpha\in\bR^{k,l}$, let
$F_{k,l}(\alpha)$ be the standard representation of
 the Heisenberg algebra with generators $j^a(n)$ and relations
$[j^a(n),j^b(m)]=n\delta_{n+m,0}\eta^{ab}id$
($a,b=1,..,k+l$, $n,m\in\bZ$,
$\eta=diag(+,..,+,-,..-)$ with $k\ +$ and $l\ -$). Here the $j^a(0)$ acts by
the scalar $\alpha^a$.
The canonical generator of the module is denoted by $|\alpha\rgl$.
We now regard each $F_{k,l}(\alpha)$ as a $Vir$-module in which $Vir$ acts by
$L(z)=\frac{1}{2}:j^a(z)j^b(z):\eta^{ab}$. This module has the standard
hermitean form \cite{FGZ}.

\be{thm}\cite{FGZ}\lb{fgz}
For any $\alpha\in\bR^{25,1}$,
the linear maps $\nu_1,\nu_2$ for the module $F_{25,1}(\alpha)$
 are bijections, and
\eq
H^{\frac{\infty}{2}+p}(Vir, F_{25,1}(\alpha))=\left\{\be{array}{ll}
\delta_{\alpha,0}\bC\bone\otimes|0\rgl & \mbox{if $p=0$}\\
\nu_1 {\frak{P}}( F_{25,1}(\alpha)) & \mbox{if $p=1$.}\\
\nu_2 {\frak{P}}( F_{25,1}(\alpha)) & \mbox{if $p=2$.}\\
\delta_{\alpha,0}\bC c(-3)c(-2)c(-1)\bone\otimes|0\rgl  & \mbox{if $p=3$}\\
0& \mbox{otherwise.}\end{array}\right.
\en
\end{thm}
For a proof, see the original reference.

Denote the highest vector of the $Vir$-modules $L(\kappa,h)$ or $M(\kappa,h)$
by $|\kappa,h\rgl$.
\be{lem}\lb{LFlemma}
For any $h\geq 0$ and $\beta\in\bR^{r-1,1}$,
the linear maps $\nu_1,\nu_2$ for the module\\
 $L(26-r,h)\otimes F_{r-1,1}(\beta)$ are bijections, and
\eq\lb{LF}
H^{\frac{\infty}{2}+p}(Vir, L(26-r,h)\otimes
F_{r-1,1}(\beta))=\left\{\be{array}{ll}
\delta_{h,0}\delta_{\beta,0}\bC\bone\otimes|26-r,0\rgl\otimes|\beta\rgl &
\mbox{if $p=0$}\\
\nu_1 {\frak{P}}( L(26-r,h)\otimes F_{r-1,1}(\beta)) & \mbox{if $p=1$.}\\
\nu_2 {\frak{P}}(L(26-r,h)\otimes F_{r-1,1}(\beta)) & \mbox{if $p=2$.}\\
\delta_{h,0}\delta_{\beta,0}
\bC c(-3)c(-2)c(-1)\bone\otimes|26-r,0\rgl\otimes|\beta\rgl
& \mbox{if $p=3$}\\
0& \mbox{otherwise.}\end{array}\right.
\en
\end{lem}
Proof: The case $r=26$ is just Theorem \ref{fgz}. So let's assume $r<26$.
The case $h=0,\beta=0$ can be easily checked by hand. So let's assume that
either $h$ or $\beta$ is nonzero. We claim that any irreducible module
$L(26-r,h)$ is direct summand in some $F_{26-r,0}(\alpha)$. Choose $\gamma$ so
that $\frac{\gamma\cdot\gamma}{2}=h$, and
we've a homomorphism $M(26-r,h)\ra F_{26-r,0}(\gamma)$ with
$|26-r,h\rgl\mapsto|\gamma\rgl$.
Since $F_{26-r,0}(\gamma)$ is unitarizable, the image must also be
unitarizable. Thus it
must also be an irreducible direct summand.

Now observe that both $H^{\frac{\infty}{2}+p}(Vir,-)$ and  ${\frak{P}}( -)$ are
exact
with respect to direct sum. Since $F_{26-r,0}(\gamma)\otimes F_{r-1,1}(\beta)$
is isomorphic to
$F_{25,1}(\alpha)$ for some $\alpha\neq0$, we see that Theorem \ref{fgz}
implies \erf{LF}. $\Box$

Proof of  Theorem \ref{general}: By assumption, $O$ is unitarizable and hence
is a direct sum of
$L(26-r,h)$, $h\geq0$. On the other hand
$V_\Lambda=\oplus_{\beta\in\Lambda}F_{r-1,1}(\beta)$
as $Vir$-modules, if we choose an identification
$\bR^{r-1,1}=\Lambda\otimes_{\bf Z}\bR$.
Now the theorem follows from Lemma \ref{LFlemma} and the fact that both
$H^{\frac{\infty}{2}+p}(Vir,-)$ and  ${\frak{P}}( -)$ are exact
with respect to direct sum. $\Box$

\subsection{Some applications}

The BV algebra $H^*(O\otimes O(V_\Lambda))$ in Theorem \ref{general2} is graded
by $\Lambda$
because as a $Vir$-module, $O(V_\Lambda)\cong V_\Lambda\cong
\oplus_{\alpha\in\Lambda}F_{1,1}(\alpha)$ is graded by $\Lambda$. In
particular, we've a decomposition of the Moonshine cohomology
${\frak{M}}^*(O(V_\Lambda))=\oplus_{\alpha\in\Lambda}{\frak{M}}^*(O(V_\Lambda))_\alpha$. Since each $F_{1,1}(\alpha)$ is tame as $Vir$-module, and since $O(V)$ is also tame, the graded dimensions $dim\ {\frak{M}}^*(O(V_\Lambda))_\alpha$ (see Theorem \ref{moonshine}) are in fact
finite. Thus the ${\frak{M}}^*(O(V_\Lambda))_\alpha$ are finite dimensional
representations of the group $F_1$. We'll compute the dimensions using our
results above together with well-known techniques in semi-infinite cohomology
theory (see \cite{FGZ}\cite{LZ2}).

Since each $O(V)\otimes F_{1,1}(\alpha)$ is hermitean, there is an induced
nondegenerate hermitean form on
the cohomology ${\frak{M}}^1(O(V_\Lambda))_\alpha$ (see below). We will compute
the signature of
this hermitean form and show that it's positive definite for nonzero $\alpha$
(``no-ghost theorem'').
We will restrict ourselves to the case when $\Lambda$ is a rank 2 hyperbolic
lattice.
How to generalize our computations to other Lorentzian lattices
will become clear, and is left as an exercise to the readers.

Since we're only interested in dimensions and signatures of cohomology, it's
enough to
work with the additive version of our results Theorem \ref{general} and Lemma
\ref{LFlemma}.
We begin by introducing one other tool: the notion of relative semi-infinite
cohomology. We
refer to readers to original references for details.

Let $C_\Delta^{\frac{\infty}{2}+*}(Vir,M)$ be the subspace of the $Vir$-module
$C^{\frac{\infty}{2}+*}(Vir,M)$ annihilated by $b(1)$ and $L_0$. Because
$[d_M,b(1)]=L_0$,
this subspace is a complex with differential $d_M$. Call this subspace the
relative complex, and
its cohomology $H_\Delta^{\frac{\infty}{2}+*}(Vir,M)$ the relative cohomology.
Note that if $M$ is tame and its weight $||\cdot||$ is bounded from below, then
 $C_\Delta^{\frac{\infty}{2}+*}(Vir,M)$
is finite dimensional.
Relative cohomology is an important tool for studying (absolute) semi-infinite
cohomology.
For example, technically to prove Theorem \ref{fgz}, we'd have to {\it first}
prove a similar vanishing
theorem for relative cohomology. In this paper, we manage to prove all our
results so far without
using it. However for computing dimensions and signatures, relative cohomology
is indispensable.

In this section, we'll be interested in
$H_\Delta^{\frac{\infty}{2}+*}(Vir,V\otimes F_{1,1}(\alpha))$,
$\alpha\in\Lambda$. For simplicity, we'll abbreviate the absolute complex
$C^{\frac{\infty}{2}+*}(Vir,V\otimes F_{1,1}(\alpha))$ simply as
$C^{\frac{\infty}{2}+*}(\alpha)$. Similar notations apply
to the absolute cohomology, the relative complex, and the relative cohomology.
\be{lem}\lb{rellemma}
\eq
H_\Delta^{\frac{\infty}{2}+p}(\alpha)\cong\left\{\be{array}{ll}
\delta_{\alpha,0}\bC & \mbox{if $p=0,2$}\\
H^{\frac{\infty}{2}+1}(\alpha) & \mbox{if $p=1$.}\\
0& \mbox{otherwise.}\end{array}\right.
\en
\end{lem}
Proof: The case $\alpha=0$ is a trivial exercise. There's a long exact sequence
(see \cite{LZ2} for details):
\eq
\cdots\ra H_\Delta^{\frac{\infty}{2}+p}(\alpha)\ra
H^{\frac{\infty}{2}+p}(\alpha)
\ra H_\Delta^{\frac{\infty}{2}+p-1}(\alpha)
\ra H_\Delta^{\frac{\infty}{2}+p+1}(\alpha)\ra\cdots.
\en
Assume $\alpha$ nonzero. By decomposing $V$ in terms of its irreducible
submodules $L(24,h)$ and applying Lemma \ref{LFlemma} (for $r=2$), we see that
the long exact sequence above degenerates
and yields the desired result. $\Box$

\be{cor}\lb{relcoho}
${\frak{M}}^1( O(V_\Lambda))_\alpha\cong
H_\Delta^{\frac{\infty}{2}+1}(\alpha)$.
\end{cor}
Thus to compute the graded dimensions of the Moonshine cohomology in degree 1
(which is
a $\Lambda$-graded Lie algebra) for the conformal QOA $O(V_\Lambda)$, it's
enough to compute
the corresponding degree 1 relative cohomology.
For a tame hermitian $||\cdot||$-graded vector space $M$, let $ch_qM=\sum dim\
M[n]q^n$,
$sign_qM=\sum sign\ M[n]q^n$.

By the Euler-Poincar\'{e} Principle, we have
\eq\lb{dim}
\sum_i (-1)^i dim\ H_\Delta^{\frac{\infty}{2}+i}(\alpha)=
\sum_i (-1)^i dim\ C_\Delta^{\frac{\infty}{2}+i}(\alpha).
\en
Note that $Ker\ b(1)=Im\ b(1)$. Thus the RHS of \erf{dim} is just the constant
term of the following
$q$-series:
\eqa
ch_qV\cdot ch_q F_{1,1}(\alpha)\cdot\sum_i(-1)^i ch_qb(1)\wedge^{i+1}
&=&q(j(q)-744)\cdot q^{\frac{\alpha\cdot\alpha}{2}} \varphi(q)^{-2} \cdot
(-1)q^{-1}\varphi(q)^2\nnb\\
&=&-q^{\frac{\alpha\cdot\alpha}{2}}(j(q)-744)
\ena
where $\varphi(q)=\Pi_{n>0}(1-q^n)$.
For $\alpha=0$, the constant term on this RHS is zero -- consistent with the
fact that
$dim\ H_\Delta^{\frac{\infty}{2}+i}(0)$, $i=0,1,2$, are respectively 1,2,1. For
$\alpha$ nonzero,
combining Lemma \ref{rellemma}, Corollary \ref{relcoho}, and eqn. \erf{dim}, we
get
\eq
dim\ {\frak{M}}^1( O(V_\Lambda))_\alpha=Res_q
q^{\frac{\alpha\cdot\alpha}{2}-1}(j(q)-744).
\en

The space $b(1)\wedge^*$ has a unique hermitean form $\lgl\cdot,\cdot\rgl$
such that $\lgl\bone,c(-3)c(-1)\bone\rgl=1$, $b(n)^\dagger=b(-n+2)$, and
$c(n)^\dagger=c(-n-4)$ \cite{FGZ}.
This makes each complex $C_\Delta^{\frac{\infty}{2}+*}(\alpha)$ a hermitian
space such that
the semi-infinite differential $d$ (dropping the subscript) is self-adjoint.
By the Euler-Poincar\'{e} Principle for signature, we have
\eq
\sum_i sign\ H_\Delta^{\frac{\infty}{2}+i}(\alpha)=
\sum_i sign\ C_\Delta^{\frac{\infty}{2}+i}(\alpha).
\en
As before, the RHS here is the constant term of the following
$q$-series:
\eqa
sign_qV\cdot sign_q F_{1,1}(\alpha)\cdot\sum_i sign_qb(1)\wedge^{i+1}
&=&q(j(q)-744)\cdot q^{\frac{\alpha\cdot\alpha}{2}} \lambda(q)^{-1} \cdot
q^{-1}\lambda(q)\nnb\\
&=&q^{\frac{\alpha\cdot\alpha}{2}}(j(q)-744)
\ena
where $\lambda(q)=\Pi_{n>0}(1-q^n)(1+q^n)$. Thus we conclude that for $\alpha$
nonzero:
\eq
sign\ H_\Delta^{\frac{\infty}{2}+1}(\alpha)=dim\
H_\Delta^{\frac{\infty}{2}+1}(\alpha)
\en
which is the statement of the no-ghost theorem. The interested readers should
compare our results with the results of Borcherds \cite{Bor2} on graded
dimensions and positive definiteness for the Monster Lie algebra.

\subsection{Remarks}

In the light of our results, a number of interesting questions come to mind. We
discuss two of
them.

1. Recall that the Virasoro element $\omega$ of the Moonshine VOA $V$ is an
$F_1$-invariant element
of $V$. This means that $Y(\omega,z)$ is an $F_1$-invariant quantum operator,
and hence
generates an $F_1$-invariant subalgebra in $O(V)$. It's easy to show that this
subalgebra
is isomorphic to $O_{24}(L)$. By unitarizability, $O_{24}(L)$ is a direct
summand in $O(V)$
as $Vir$-module. Thus for any conformal QOA $O$ of central charge 2, the
Moonshine cohomology
${\frak{M}}^*(O)$, as a BV algebra, has a canonical subalgebra
$H^*(O_{24}(L)\otimes O)$. In
particular, the vanishing theorem \ref{moonvanish} holds for this algebra.
Also, $H^2(O_{24}(L)\otimes O)$ is the adjoint module over
the Lie algebra $H^1(O_{24}(L)\otimes O)$.

{\bf Problem:} For any rank two hyperbolic lattice $\Lambda$, study the
$\Lambda$-graded Lie subalgebra\\ $H^1(O_{24}(L)\otimes O(V_\Lambda))$ of
${\frak{M}}^1(O(V_\Lambda))$.

2. Recall that the Monster finite group $F_1$ acts naturally on
${\frak{M}}^*(O)$; thus, $F_1$ acts as a group
of invertible natural transformations of ${\frak{M}}^*$ as a functor. Since
the full automorphism group of the Moonshine VOA $V$, hence of the conformal
QOA $O(V)$, is isomorphic to $F_1$, the following seems quite plausible:

{\bf Conjecture: } The group of natural automorphisms of the Moonshine
cohomology functor is isomorphic to $F_1$.

\noi{\footnotesize DEPARTMENT OF MATHEMATICS, HARVARD UNIVERSITY
CAMBRIDGE, MA 02138.\\ lian$@$math.harvard.edu}\\
\noi{\footnotesize DEPARTMENT OF MATHEMATICS, YALE UNIVERSITY
NEW HAVEN, CT 06520.\\ gregg$@$math.yale.edu}

\end{document}